%Paper: alg-geom/9510015
%From: "Steven L. Kleiman" <kleiman@math.mit.edu>
%Date: Mon, 30 Oct 1995 14:47:58 -0500 (EST)
%Date (revised): Mon, 30 Oct 1995 15:13:35 -0500 (EST)
%Date (revised): Mon, 30 Oct 1995 15:22:20 -0500 (EST)
%Date (revised): Fri, 1 Mar 1996 07:50:22 -0500 (EST)
%Date (revised): Fri, 1 Mar 1996 07:56:22 -0500 (EST)

  %%%%%%%%%%%%%%%%%%%%%%%%%%%%%%%%%%%%%%%%%%%%%%%%%%%%%%%%%%%%%%%%%%%%%%
  %                                                                    %
  % Rational curves of degree at most 9 on a general quintic threefold %
  %                                                                    %
  %               Trygve JOHNSEN and Steven L. KLEIMAN                 %
  %                                                                    %
  %                          alg-geom/9510015
  %                                                                    %
  %%%%%%%%%%%%%%%%%%%%%%%%%%%%%%%%%%%%%%%%%%%%%%%%%%%%%%%%%%%%%%%%%%%%%%
%%
%% Plain Tex
%%
%% SOME CHOICES
\def\TheMagstep{\magstep1}	% Normal magnification
		% Changed to \magstep0 by \DoublepageOutput{TRUE}
\def\PaperSize{letter}		% \PaperSize is used to
% \def\PaperSize{AFour}		%  center text on page
%%%%%
%%% To get two pages side by side in landscape mode on an ordinary
%%    sheet via dvips and  PostScript, enable the next command.
% \def\DoublepageOutput{TRUE}
%%	END of CHOICES

\def\TheABSTRACT{
 We prove the following form of the Clemens conjecture in low degree.
Let $d\le9$, and let $F$ be a general quintic threefold in $\IP^4$.
Then (1)~the Hilbert scheme of rational, smooth and irreducible curves
of degree $d$ on $F$ is finite, nonempty, and reduced; moreover, each
curve is embedded in $F$ with normal bundle $\O(-1)\oplus\O(-1)$, and
in $\IP^4$ with maximal rank.  (2)~On $F$, there are no rational,
singular, reduced and irreducible curves of degree $d$, except for the
17,601,000 six-nodal plane quintics (found by Vainsencher).  (3)~On
$F$, there are no connected, reduced and reducible curves of degree $d$
with rational components.}

%% Headlines
 \def\firstheadline{\hfil}
 \def\oddheadline{\rm
  \rlap{January 24, 1996}
  \hfil\shorttitle\hfil\llap{\the\pageno}}
 \def\evenheadline{\rm\rlap{\number\pageno}\hfil
 	\author\hfil \llap{alg-geom/9510015}}

\let\:=\colon \def\IP{{\bf P}} \let\x=\times \let\ox=\otimes
\def\O{{\cal O}} \def\I{{\cal I}} \def\A{{\cal A}} 
\def\M{{\cal M}}
\let\To=\longrightarrow \def\TO#1{\buildrel#1\over\To}
\def\and{\hbox{ and }}

\def\mathopdef#1{\expandafter\def
 \csname#1\endcsname{\mathop{\rm #1}\nolimits}}
\def\NoOp{*!*}
\def\NextOp#1 {\def\TheOp{#1}\ifx\TheOp\NoOp\let\next\relax
  \else\mathopdef{#1}\let\next\NextOp \fi \next}
\NextOp
 cod rank grad mod Ker
 *!*

			       %% FORMAT
%% Page layout
\parskip=0pt plus 1.75pt
\parindent10pt
\hsize26pc
% \vsize42pc
%\vsize=42\normalbaselineskip %42 lines per page
% -SLK-95OCT21=>
\vsize=43\normalbaselineskip %43 lines per page

	% to shorten the space above and below displays
\abovedisplayskip6pt plus6pt minus0.25pt
\belowdisplayskip6pt plus6pt minus0.25pt

\def\TRUE{TRUE}
\ifx\DoublepageOutput\TRUE \def\TheMagstep{\magstep0} \fi
\mag=\TheMagstep

% Center text on page
	% additional vertical adjustment
\newskip\vadjustskip \vadjustskip=0.5\normalbaselineskip
\def\centertext
 {\hoffset=\pgwidth \advance\hoffset-\hsize
  \advance\hoffset-2truein \divide\hoffset by 2\relax
  \voffset=\pgheight \advance\voffset-\vsize
  \advance\voffset-2truein \divide\voffset by 2\relax
  \advance\voffset\vadjustskip
 }
\newdimen\pgwidth\newdimen\pgheight
\def\letter{letter}\def\AFour{AFour}
\ifx\PaperSize\letter
 \pgwidth=8.5truein \pgheight=11truein
 \message{- Got a paper size of letter.  }\centertext
\fi
\ifx\PaperSize\AFour
 \pgwidth=210truemm \pgheight=297truemm
 \message{- Got a paper size of AFour.  }\centertext
\fi

%% Two-column landscape format
% Modified from the TeX book, p. 257.
 \newdimen\fullhsize \newbox\leftcolumn
 \def\fulline{\hbox to \fullhsize}
\def\doublepageoutput
{\let\lr=L
 \output={\if L\lr
          \global\setbox\leftcolumn=\columnbox \global\let\lr=R%
        \else \doubleformat \global\let\lr=L\fi
        \ifnum\outputpenalty>-20000 \else\dosupereject\fi}%
 \def\doubleformat{\shipout\vbox{%
        \fulline{\hfil\hfil\box\leftcolumn\hfil\columnbox\hfil\hfil}%
}%
 }%
 \def\columnbox{\vbox
   {\makeheadline\pagebody\makefootline\advancepageno}%
   }
 \fullhsize=\pgheight \hoffset=-1truein
 \voffset=\pgwidth \advance\voffset-\vsize
  \advance\voffset-2truein \divide\voffset by 2
  \advance\voffset\vadjustskip
\let\firstheadline=\hfil

\null\vfill\nopagenumbers\eject\pageno=1\relax % to put page on right
}
%\ifx\DoublepageOutput\TRUE \doublepageoutput \fi

\catcode`\@=11
%% Modification of the PLAIN footnote macro for 8pt
\def\vfootnote#1{\insert\footins\bgroup
 \eightpoint %% only change
 \interlinepenalty\interfootnotelinepenalty
  \splittopskip\ht\strutbox % top baseline for broken footnotes
  \splitmaxdepth\dp\strutbox \floatingpenalty\@MM
  \leftskip\z@skip \rightskip\z@skip \spaceskip\z@skip \xspaceskip\z@skip
  \textindent{#1}\footstrut\futurelet\next\fo@t}

%% Headlines %%
 % Switchs for dates in the headlines
\newif\ifdates	\newif\ifdateonpageone
% From the TeX book p. 406
\def\today{\ifcase\month\or
 January\or February\or March\or April\or May\or June\or
 July\or August\or September\or October\or November\or December\fi
 \space\number\day, \number\year}

\nopagenumbers
\headline={%
 \eightpoint
  \ifnum\pageno=1\firstheadline
  \else
    \ifodd\pageno\oddheadline
    \else\evenheadline\fi
  \fi
}

%% To start a new section.
%%  Checks that there's room enough on the page.
 \newcount\sectno \sectno=0
 \newskip\sectskipamount \sectskipamount=0pt plus30pt
 \def\newsect #1\par{\displayno=0 %\proclaimno=0
   \advance\sectno by 1
   \vskip\sectskipamount\penalty-250\vskip-\sectskipamount
   \bigskip	%	\smallskip \medskip
   \centerline{\bf \number\sectno. #1}\nobreak
   \medskip	%	\smallskip
   \message{#1 }
}
			      %% STYLE %%
\hyphenation{par-a-met-ri-za-tion par-a-met-ri-zed}
%% For Roman paranthetical material in nonRoman text
\def\(#1){{\rm(#1)}}\let\leftp=(
\def\activeleftp{\catcode`\(=\active}
{\activeleftp\gdef({\ifmmode\let\next=\leftp \else\let\next=\(\fi\next}}

%% SECTIONING
%% Articles
\def\artkey #1 {{\bf (\number\sectno.#1)}\enspace}

 \let\rmk=\rem

\def\proclaim#1 #2 {\medbreak{\bf#1 }\artkey#2 \bgroup\it\activeleftp}
\def\endproclaim{\egroup\medbreak}
\def\pf{\endproclaim{\bf Proof.}\enspace}
\def\prop{\proclaim Proposition } 
\def\lem{\proclaim Lemma }
\def\thm{\proclaim Theorem }
\def\cor{\proclaim Corollary }

%% Display numbers
%% Convient macro for display numbers
\newcount\displayno
\def\eqlt#1${\global\advance\displayno by 1
 \expandafter\xdef
  \csname \the\sectno#1\endcsname{\number\displayno}
 \eqno\hbox{\rm(\the\sectno-\number\displayno)}$}
\def\tjitem #1 {\global\advance\displayno by 1
 \expandafter\xdef
  \csname \the\sectno#1\endcsname{\number\displayno}
 \item{\the\sectno-\number\displayno} }
\def\eqalignlt#1\cr{\global\advance\displayno by 1
 \expandafter\xdef
  \csname \the\sectno#1\endcsname{\number\displayno}
 \eqno\hbox{\rm(\the\sectno-\number\displayno)}\cr}
\def\disp#1#2{(#1-\csname#1#2\endcsname)}
\def\Cs#1){\unskip~{\rm(\number\sectno.#1)}}
\def\Cd#1){\unskip~{\rm(\number\sectno-\csname\the\sectno#1\endcsname)}}
\def\Co#1-#2){\unskip~{\rm(#1-\csname#1#2\endcsname)}}

%% Redefine \item to give the indentation of AMSTeX and the roman font
\def\item#1 {\par\indent\indent\indent \hangindent3\parindent
 \llap{\rm (#1)\enspace}\ignorespaces}
%% Define a similar macro without the hanging indentation for assertions
%% and  that starts each part with an ordinary \parindent
 \def\part#1 {\par{\rm (#1)\enspace}\ignorespaces}

%% REFERENCING %%
%%  Macros for introducing the reference keys in order
 \newif\ifproofing \proofingfalse % Decide between Alpha and Num keys.
 \newcount\refno	 \refno=0
 \def\MakeKey{\advance\refno by 1 \expandafter\xdef
 	\csname\TheKey\endcsname{{%
	\ifproofing\TheKey\else\number\refno\fi}}\NextKey}
 \def\NextKey#1 {\def\TheKey{#1}\ifx\TheKey\NoKey\let\next\relax
  \else\let\next\MakeKey \fi \next}
 \def\NoKey{*!*}
 \def\RefKeys #1\endRefKeys{\expandafter\NextKey #1 *!* }

%% \cite code from AMSTeX, modified so first argument is a control sequence
% First redefine ties to give a thinspace
\def\UThin{\penalty\@M \thinspace\ignorespaces}
\def\relaxnext@{\let\next\relax}
\def\cite#1{\relaxnext@
 \def\nextiii@##1,##2\end@{\unskip\space{\rm[\SetKey{##1},\let~=\UThin##2]}}%
 \in@,{#1}\ifin@\def\next{\nextiii@#1\end@}\else
 \def\next{\unskip\space{\rm[\SetKey{#1}]}}\fi\next}
\newif\ifin@
\def\in@#1#2{\def\in@@##1#1##2##3\in@@
 {\ifx\in@##2\in@false\else\in@true\fi}%
 \in@@#2#1\in@\in@@} \def\SetKey#1{{\bf\csname#1\endcsname}}

\let\texttilde=\~
\def\~{\ifmmode\let\next=\widetilde\else\let\next=\texttilde\fi\next}

%%  ``Ties'' with a \thinspace for page numbers
\def\p.{\unskip\space p.\UThin}
\def\pp.{\unskip\space pp.\UThin}
%% Roman parentheses for use within text set in other fonts
%% taken from the TeXbook, p.409, but modified to use the italic correction
%\def\(#1){{\rm(}#1\/{\rm)}}

%% FONTS
 \font\twelvebf=cmbx12
 \font\smc=cmcsc10

%% Eight point type for footnotes and references
\def\eightpoint{\eightpointfonts
 \setbox\strutbox\hbox{\vrule height7\p@ depth2\p@ width\z@}%
 \eightpointparameters\eightpointfamilies
 \normalbaselines\rm
 }
\def\eightpointparameters{%
 \normalbaselineskip9\p@
 \abovedisplayskip9\p@ plus2.4\p@ minus6.2\p@
 \belowdisplayskip9\p@ plus2.4\p@ minus6.2\p@
 \abovedisplayshortskip\z@ plus2.4\p@
 \belowdisplayshortskip5.6\p@ plus2.4\p@ minus3.2\p@
 \skewchar\eighti='177 \skewchar\sixi='177
 \skewchar\eightsy='60 \skewchar\sixsy='60
 \hyphenchar\eighttt=-1
 }
\newfam\smcfam
\def\eightpointfonts{%
 \font\eightrm=cmr8 \font\sixrm=cmr6
 \font\eightbf=cmbx8 \font\sixbf=cmbx6
 \font\eightit=cmti8
% \font\eightsmc=cmcsc10 at 10truept
\font\eightsmc=cmcsc8
 \font\eighti=cmmi8 \font\sixi=cmmi6
 \font\eightsy=cmsy8 \font\sixsy=cmsy6
 \font\eightsl=cmsl8 \font\eighttt=cmtt8}
\def\eightpointfamilies{%
 \textfont\z@\eightrm \scriptfont\z@\sixrm  \scriptscriptfont\z@\fiverm
 \textfont\@ne\eighti \scriptfont\@ne\sixi  \scriptscriptfont\@ne\fivei
 \textfont\tw@\eightsy \scriptfont\tw@\sixsy \scriptscriptfont\tw@\fivesy
 \textfont\thr@@\tenex \scriptfont\thr@@\tenex\scriptscriptfont\thr@@\tenex
 \textfont\itfam\eightit	\def\it{\fam\itfam\eightit}%
 \textfont\slfam\eightsl	\def\sl{\fam\slfam\eightsl}%
 \textfont\ttfam\eighttt	\def\tt{\fam\ttfam\eighttt}%
 \textfont\smcfam\eightsmc	\def\smc{\fam\smcfam\eightsmc}%
 \textfont\bffam\eightbf \scriptfont\bffam\sixbf
   \scriptscriptfont\bffam\fivebf	\def\bf{\fam\bffam\eightbf}%
 \def\rm{\fam0\eightrm}%
% \tt \ttglue=0.5em plus0.25em minus0.15em
 }

\catcode`\@=13

\ifx\DoublepageOutput\TRUE \doublepageoutput \fi

%% Macros for setting references
 \def\SetRef#1 #2\par{%
   \hang\llap{[\csname#1\endcsname]\enspace}%
   \ignorespaces#2\unskip.\endgraf}
 \newbox\keybox \setbox\keybox=\hbox{[18]\enspace}
 \newdimen\keyindent \keyindent=\wd\keybox
 \def\references{\vskip-\smallskipamount
  \bgroup   \eightpoint   \frenchspacing
   \parindent=\keyindent  \parskip=\smallskipamount
   \everypar={\SetRef}}
 \def\endreferences{\egroup}
 \def\nocomma{\def\eatcomma##1,{}\expandafter\eatcomma}

%% Papers and books
  \def\paper{\unskip, \bgroup\it}
 \def\paperinfo{\unskip, \egroup}
 \def\preprint{\unskip, \egroup Preprint }
 \def\inbook#1\bookinfo#2\publ#3\yr#4\pages#5
   {\unskip, \egroup in ``#1\unskip,'' #2\unskip, #3\unskip,
   #4\unskip, pp.~#5}
 \def\at.{.\spacefactor3000}
 \def\book{\unskip, ``}

 \def\bookinfo#1{\unskip," #1}

\def\MirrorMan{\unskip, \egroup in ``Essays on mirror manifolds,''
S.-T.  Yau (Ed.), International Press, Hong Kong (1992), }
 \def\Tromso{\unskip, \egroup in ``Algebraic Geometry $\bullet$
Proceedings, Troms\o, Norway 1977," L. Olsen (Ed.), Springer Lecture
Notes in Math. {\bf 687} (1978), }
 \def\RdiPapa{\unskip, \egroup in ``Space Curves $\bullet$ Proceedings,
Rocca di Pappa 1985," F. Ghione, C. Peskine, E. Sernesi (Eds.),
Springer Lecture Notes in Math. {\bf 1266} (1986), }
 \def\Zeuthen{\unskip, \egroup in ``Enumerative Algebraic
Geometry$\bullet$ Proceedings of the 1989 Zeuthen Symposium," S.
L. Kleiman and A. Thorup (Eds.), CONM {\bf 123} (1991), }

%% SERIALS
 \def\serial#1#2{\expandafter\def\csname#1\endcsname ##1 ##2 ##3
  {\unskip, \egroup #2 {\bf##1} (##2), ##3}}
 \serial{acta}{Acta Math.}
 \serial{crasp}{C. R. Acad. Sci. Paris}
 \serial{comp}{Comp. Math.}
 \serial{duke}{Duke Math. J.}
 \serial{pmihes}{Publ. Math. I.H.E.S.}
 \serial{invent}{Invent. Math.}
 \serial{jag}{J. Alg. Geom.}
 \serial{crelle}{J. reine ang. Math.}
 \serial{slnm}{Springer Lecture Notes in Math}
 \serial{slnp}{Springer Lecture Notes in Physics}
 \serial{ma}{Math. Ann.}
 \serial{mm}{Manuscripta Math.}
 \serial{ms}{Math. Scand.}
 \serial{mathz}{Math. Z.}
 \serial{np}{Nuclear Physics. B}
 \serial{qjm}{Quart. J. Math. Oxford Ser. (2)}

 \datestrue			%% want dates in headlines

 \RefKeys
 ACGH BE Cetal Cil C C86 CJ dA ES GP GLP Gu H Hirsch Huyb Katz Kont L
LY LT N Oguiso Ramella Rath V Verdier Westhoff
 \endRefKeys

%%% DOCUMENT
%% Top Matter with the format of the AMS CNM.STY
\def\author{T. Johnsen and S. Kleiman}
\def\shorttitle{Rational curves of degree at most 9}

\long\def\topstuff{
\null\vskip24pt plus 12pt minus 12pt
\twelvebf			%% 12pt bold for title
\centerline{Rational curves of degree at most 9}\medskip
\centerline{on a general quintic threefold}

 \rm				%% return to 10pt Roman

%% Subject classification and acknowledgements
 \footnote{}{\noindent 1980 {\it Mathematics Subject Classification}
   (1985 {\it Revision}).  Primary 14J30; Secondary 14H45, 14N10.}
 \footnote{}{It is a pleasure to thank Stein Arild Str\o mme for
discussions about the higher-order deformation theory involved in the
proof of Lemma~(2.4), and to thank Kristian Ranestad for suggesting the
use of Hirschowitz's ``m\'ethode d'Horace'' (see \cite{Hirsch}) to prove
an earlier version of Lemma~(4.5).}

%% Authors' names and support
 \vskip15pt plus6pt minus3pt
 \centerline{\smc
 Trygve JOHNSEN\footnote{$^{1}$}{%
    Supported in part by the Norwegian Research Council for Science and
the Humanities.  This author is grateful to the MIT Department of
Mathematics for its hospitality from January 1 to July 31, 1993, when
this work was started.}
 and Steven L.~KLEIMAN\footnote{$^{2}$}{%
    Supported in part by NSF grants 9106444-DMS and 9400918-DMS.  This
author is grateful to the Institute of Advanced Studies of the
Norwegian Academy of Sciences for its hospitality from January 31 to
February 25, 1994 and to the Mathematical Institute of the University
of Copenhagen for its hospitality during the summer of 1995; during
these periods this work was continued and nearly completed.}}

%% Preprint identification
%\vskip12pt plus12pt minus12pt
% \centerline{alg-geom/9510015}

%% Abstract
\vskip15pt plus12pt minus12pt
{\parindent=24pt \narrower \noindent \eightpoint
 {\smc Abstract.}\enspace\TheABSTRACT \par }
 %\vskip6pt plus12pt minus6pt
 } %\end of topstuff

 \topstuff

\newsect	Introduction

Around ten years ago, Clemens posed a conjecture about the rational
curves on a general quintic threefold $F$ in complex $\IP^4$.
Immediately afterwards, S.~Katz \cite{Katz} considered the conjecture in
the following form: {\it the Hilbert scheme of rational, smooth and
irreducible curves of degree $d$ on $F$ is finite, nonempty and
reduced; in fact, each curve is embedded with normal bundle
$\O(-1)\oplus\O(-1)$}.  Katz proved this statement for $d\le7$.  Then
Clemens \cite{C86, \p.639} strengthened the conjecture by adding this
assertion: {\it all the rational, reduced and irreducible curves on $F$
are smooth and mutually disjoint}.  This additional assertion is not
completely true; Vainsencher \cite{V, \S6.2} proved that $F$ contains
17,601,000 six-nodal quintic plane curves.  By developing Katz's
argument through an extensive case-by-case analysis, we'll prove a
corrected form of Clemens' stronger conjecture for $d\le9$.  In
Section~2, we'll prove Katz's form of the conjecture.  In Section~3,
we'll prove that $F$ contains no singular, rational, reduced and
irreducible curves of degree $d$, except for the six-nodal plane
quintics.  Finally, in Section~4, we'll prove that $F$ contains no
connected, reduced and reducible curves of degree $d$ whose components
are all rational.

Consider the incidence scheme $I_d$ that parametrizes the pairs $(C,F)$
where $C$ is a rational, smooth and irreducible curve of degree $d$ in
$\IP^4$ and where $F$ is a quintic threefold containing $C$.  For every
$d$, Katz \cite{Katz, \p.152} proved that his form of Clemens'
conjecture is true if $I_d$ is irreducible and if there exists at least
one pair $(C,F)\in I_d$ such that $F$ is smooth and such that $C$,
viewed as $\IP^1$, is embedded with normal bundle $\O(-1)\oplus\O(-1)$.
He \cite{Katz, \p.153} proved the existence of such a pair for every
$d$ (Clemens \cite{C, Thm.~(1.27), \p.26} had just proved it for
infinitely many $d$, and Katz observed that the general case follows
via Clemens' deformation-theoretic argument from an existence result of
Mori's).  And Katz proved the irreducibility of $I_d$ for $d\le7$ by
making the following two observations: (1)~the space $M_d$ of all $C$
is irreducible for every $d$; and (2)~every fiber of the first projection
$I_d\to M_d$ is, thanks to \cite{GLP}, a projective space of the same
dimension for $d\le7$.  For $d\ge8$, however, some fibers are bigger
than others.  Nevertheless, we'll prove in Section~2 that there aren't
too many big fibers for $d=8,9$.  Thus we'll obtain the irreducibility
and so Katz's form of Clemens' conjecture; see Proposition~(2.2) and
Theorem~(2.1).  As a bonus, we'll obtain Corollary~(2.5), which asserts
that each $C$ on a general $F$ has some additional important
properties; notably, each $C$ is embedded in $\IP^4$ with maximal rank,
and the restricted twisted sheaf of differentials
$\Omega^1_{\IP^4}(1)|C$ has a certain splitting type.

Katz's form of Clemens' conjecture was recently considered
independently by Nijsse \cite{N}, by Huybrechts \cite{Huyb}, and by
Westhoff \cite{Westhoff, \p.40}.  Nijsse too proved the irreducibility
of $I_d$ for $d\le9$ by following Katz's approach, but he handled the
curves with big fibers in a rather different manner than we'll do in
Section~2.  Huybrechts and Westhoff did not follow Katz's approach, but
used more deformation theory.  Each assumed $d\le7$, and proved this
slightly weaker statement: {\it if the Hilbert scheme of rational,
smooth and irreducible curves of degree $d$ on a general $F$ is
nonempty, then it's finite and reduced}.  Also, each didn't recover the
irreducibility of $I_d$.  However, for the other Calabi--Yau threefolds
$F$ that are general complete intersections, each was able to prove the
preceding statement for correspondingly small $d$, namely, for $d\le5$
for $F$ of type (3,3), and for $d\le4$ for types (2,4), (2,2,3) and
(2,2,2,2).  In this connection, Oguiso \cite{Oguiso} proved, by
modifying Clemens' argument for quintics, that the Hilbert scheme is,
in fact, nonempty for an $F$ of type (2,4).

The above ranges of $d$, established independently by Huybrechts and by
Westhoff, can, almost certainly, be extended somewhat by proceeding
instead along the lines in Section~2.  The authors have checked the key
details, and believe the following ranges come out: $d\le7$ for types
(3,3) and (2,4), and $d\le6$ for types (2,2,3) and (2,2,2,2).  In fact,
except for the case $d=6$ and $F$ of type (2,2,2,2), the incidence
scheme $I_d$ of pairs $(C,F)$ is, almost certainly, irreducible,
generically reduced, and of the same dimension as the space $\bf P$ of
$F$.  In the exceptional case, the authors cannot rule out the presence
of a second component, whose $C$\/s lie on quadric surfaces.  In any
event, if the second projection $I_d\to{\bf P}$ is surjective, then by
Sard's theorem it is \'etale over a dense open set of $F$s, and so for
these $F$ the Hilbert scheme is finite and reduced.

To prove that a general quintic threefold does not contain curves
$C$ of a given sort, we'll form the incidence scheme of all pairs
$(C,F)$ where $F$ is an arbitrary quintic, and we'll show that this
scheme does not dominate the $\IP^{125}$ of all quintics $F$.  In
Section~3, we'll treat $C$s that are reduced and irreducible, rational
and singular.  Here, for convenience, we'll work with parametrized
mappings of $\IP^1$ into $\IP^4$, not just their images $C$.  Since
each $C$ is the image of every mapping in a four parameter family,
we'll need to show that the new incidence scheme is of dimension at
most 128.  We'll succeed in doing so for all such $C$ of degree at most
9, with one exception: the six-nodal quintic plane curves.  However,
Vainsencher \cite{V, \S6.2} proved that a general quintic threefold
contains 17,601,000 such plane curves.  Each is the intersection of the
threefold with a sixfold tangent plane.  Therefore, we're led to ask if
a general quintic threefold also contains sixteen-nodal curves of
degree 10, arising from tangent quadric surfaces.  The answer is no,
see Remark~(3.6).

Consider a connected, reduced and reducible curve $C$ of degree at most
9 with rational components.  In Section~4, we'll prove that $C$ cannot
lie on a general quintic threefold.  To do so, we'll have a
considerable advantage, gained from our work in Sections 2 and 3;
namely, we'll be able to assume that $C$ has two components, each of
which is one of the following: a line, a conic, a six-nodal plane
quintic, a twisted cubic, or a smooth curve that spans $\IP^4$ and on
which the restriction of $\Omega^1_{\IP^4}(1)$ has a certain splitting
type.  We'll parametrize the various $C$ of this sort by parametrizing
their components individually using the Hilbert scheme of $\IP^4$.
We'll bound the dimension of the incidence scheme in two steps.  First
we'll consider those $C$ whose two components meet in a subscheme of
length $n$, and prove that they form a subset of codimension at least
$n$.  Second, we'll bound the number of independent quintic threefolds
containing a given $C$.  To do so, we'll need to work harder than ever
to prove that $C$ is 6-regular.  Curiously, the regularity of $C$ seems
to depend only on the geometry of its components individually, and not
on how they intersect.

Clemens' full conjecture \cite{C86, \p.639} has one additional
assertion: the number $n_d$ of rational curves of degree $d$, on the
general quintic threefold, is divisible by $5^3\cdot d$.  This
assertion holds for $d=1,2$; indeed, $n_1$ was found by Schubert, and
$n_2$ by Katz \cite{Katz, Thm.~3.1, \p.154}.  However, Clemens hedged
on the factor of $d$, and in fact, $n_3$ is not divisible by $3$;
indeed, $n_3$ was found by Ellingsrud and Str\o mme \cite{ES} and by
Candelas, De la Ossa, Green, and Parkes \cite{Cetal}.  The latter four
were inspired by considerations of mathematical physics.  They
predicted values of $n_d$ for $d\le10$, and these values are all
divisible by $5^3$.  Moreover, they gave an algorithm yielding numbers
for all $d$.  These numbers were shown, by Lian and Yau \cite{LY}, to
be divisible by $5^3$, at least if $d$ isn't a multiple of $5$.  On the
other hand, Kontsevich \cite{Kont, \S2.2} gave a somewhat different
algorithm; although he too was inspired by mathematical physics, his
treatment is more algebraic-geometric, and its numbers clearly count
both smooth and nodal curves.  For the other four types of Calabi--Yau
complete intersections of dimension 3, Libgober and Teitelbaum
\cite{LT} followed in the footsteps of Candelas et al., and developed
algorithms.  At any rate, all this enumeration has given rise to
renewed interest in Clemens' conjecture in its full form.

Finally, a few words are in order concerning the choice of the ground
field and the meaning of the word ``general.''  As usual, we'll say
that a {\bf general} quintic threefold in $\IP^4$ has a given property if
there exists a dense Zariski open subset of $\IP^{125}$ whose points
represent threefolds with this property.  We'll consider properties
that concern curves of a fixed degree, never curves of infinitely many
different degrees simultaneously.  When the latter is done, it may be
necessary to consider the intersection of a countably infinite
collection of dense Zariski open subsets.  Such an intersection is
nonempty if the ground field is uncountable~--- for example, the
complex numbers~--- and then it is common to call a threefold
``generic'' if it's represented by a point of the intersection.

The field of complex numbers is virtually always the ground field of
choice in treatments of Clemens' conjecture.  In this paper, however,
all our work is algebraic-geometric in nature, and makes sense,
although it's not always valid, in any characteristic.  The need for
characteristic zero enters into our work in three ways: (i) our use, at
the beginning of Section~2, of Katz's theorem on the existence of
rational smooth curves on a general quintic threefold, (ii) our use of
the principle of uniform position, or the trisecant lemma, for a
singular curve to prove (3-8) via (2-3), and (iii) our use of the
Castelnuovo--Halphen bounds on the genera of singular curves, made in
the two paragraphs before the statement of Lemma~(3.2) and again to
prove (3-8).  Thus, in positive characteristic, Proposition~(2.2) is
valid, but possibly the projection $I_d\to\IP^{125}$ is not surjective
or not separable.  So a general quintic threefold contains at most
finitely many rational, smooth and irreducible curves of degree at most
9, but possibly there are none, or if there are some, then they may
move infinitesimally.  Also, the threefold may possibly contain
singular, rational, reduced and irreducible curves that are ``strange''
in the sense that every tangent line is a trisecant.  Therefore, we'll
assume implicitly from now on that the ground field is algebraically
closed and of characteristic zero.

\newsect	Smooth curves

In this section, we'll prove the following theorem, which affirms
Katz's form of Clemens' conjecture in degree $d$ at most 9.  To do so,
we'll prove Proposition~(2.2), which asserts the irreducibility of the
incidence scheme of pairs $(C,F)$ where $C$ is a rational and smooth,
reduced and irreducible curve of degree $d$ in $\IP^4$ and where $F$ is
a quintic threefold containing $C$.  As a bonus, we'll obtain
Corollary~(2.5), which gives some additional properties possessed by
the $C$ on a general $F$.  In the remaining two sections, we will
complete our treatment of Clemens' conjecture by proving two theorems
about singular curves.

 \thm1 Let $d\le9$, and let $F$ be a general quintic threefold in
$\IP^4$.  In the Hilbert scheme of $F$, form the open subscheme of
rational, smooth and irreducible curves $C$ of degree $d$.  Then this
subscheme is finite, nonempty, and reduced; in fact, each $C$ is
embedded in $F$ with normal bundle $\O_{\IP^1}(-1)\oplus\O_{\IP^1}(-1)$.
 \endproclaim

To prove the theorem, we will use the following notation: \smallbreak
 \item a $M_d$ will denote the open subscheme of the Hilbert scheme of
$\IP^4$ parametrizing the rational, smooth and irreducible curves $C$
of degree $d$;
 \item b $\IP^{125}$ will denote the projective space
parametrizing the quintic threefolds $F$ in $\IP^4$;
 \item c $I_d$ will denote the ``incidence'' subcheme of $M_d\times
\IP^{125}$ of pairs $(C,F)$ such that $C\subset F$;
 \item d  $M_{d,i}$ will denote the locally closed subset of  $M_d$
parametrizing the curves $C$ such that $h^1(\I_C(5))=i$ where $\I_C$
denotes the ideal of $C$ in $\IP^4$;
 \item e  $I_{d,i}$ will denote the preimage in $I_d$ of $M_{d,i}$.
 \smallbreak\noindent
 By virtue of Katz's work (see the beginning of Section~1),
Theorem~\Cs1) is a consequence of the following proposition, whose
proof will occupy most of the rest of this section.

\prop2 For $d\le9$, the incidence scheme $I_d$ is irreducible of
dimension $125$, and it projects onto $M_d$.
 \endproclaim

First, let $d$ be arbitrary.  Recall from \cite{C, \pp.25--6} and
\cite{Katz, \p.152} that $M_d$ is irreducible of dimension $5d+1$ and that
every component of $I_d$ has dimension at least 125.  On the other
hand, consider $I_{d,i}$.  Its fiber over a curve $C$ in $M_{d,i}$ is a
projective space of dimension $h^0(\I_C(5))-1$.  This dimension can be
computed using the standard exact sequence,
	$$0\to H^0(\I_C(5))\to H^0(\O_{\IP^4}(5))\to
	H^0(\O_C(5))\to H^1(\I_C(5))\to 0.$$
 Therefore,
	$$\dim I_{d,i} = \dim M_{d,i} + (126-(5d+1)+i) -1.\eqlt1b$$
 If $d\le7$, then the theorem on \p.492 of \cite{GLP} yields
$M_{d,0}=M_d$, and it follows that $I_d$ is irreducible of dimension
$125$.  For $d=8,9$ the following lemma implies that $M_{d,0}$ is open
and nonempty, that $\dim I_{d,0} = 125$, and that $\dim I_{d,i}<125$ for
$i>0$; hence, $I_{d,0}$ is irreducible of dimension $125$, and its
closure is $I_d$.  Thus Proposition~(2.2) is proved, given Lemma~\Cs3).

\lem3 For $d=8,9$, if $i>0$ and if $M_{d,i}$ is nonempty, then
	 $$\cod(M_{d,i},M_d)>i.$$
 \endproclaim

To prove the lemma, we begin with some general observations.  Fix an
arbitrary rational, smooth and irreducible curve $C$ of any degree $d$
in $\IP^4$.  First, if $d\ge3$, then $C$ can lie in no plane because its
arithmetic genus $p_aC$ vanishes.  Moreover, if $d\ge4$, then $C$ can lie
in no 2-dimensional quadric cone by \cite{H, Exr.~2.9, \p.384}.
Second, if $C$ lies in a hyperplane $G$, say with ideal $\I_{C/G}$, then
	$$h^1(\I_C(5))=h^1(\I_{C/G}(5)).\eqlt1a$$
 Indeed, this equality results from the exact sequence of twisted
ideals,
	 $$0\to\I_G(5)\to\I_C(5)\to\I_{C/G}(5)\to0,$$
 because its first term is equal to $\O_{\IP^4}(4)$.

Assume $d=8$.  By the theorem on \p.492 of \cite{GLP}, if $C$ does not
lie in a hyperplane, then $C$ is 6-regular, and so $h^1(\I_C(5))=0$.
Hence, $M_8-M_{8,0}$ is contained in the closed set $N_8$ of curves in
hyperplanes.  Clearly,
$$\cod(N_8,M_8)=(5\times8+1)-(4\times8+4)=8-3=5.$$ In particular,
$\cod(M_{8,1},M_8)>1$ as asserted.

Suppose that $h^1(\I_C(5))\ge2$.  Say that $C$ lies in the hyperplane
$G$. Then $h^1(\I_{C/G}(5))\ge2$ by \Cd1a).  Hence, the table on \p.504
of \cite{GLP} says that $h^1(\I_{C/G}(5))=5$.  Therefore, if $M_{8,i}$
is nonempty, then $i$ is $0$, $1$, or $5$.  Now, the table also says
that $C$ lies in a smooth quadric surface $Q$ contained in $G$.  Since
$C$ is smooth, rational, and of degree 8, it is of type (7,1) on $Q$.
Hence, $C$ varies in a system of dimension $8\times2-1$, or 15, on $Q$.
Obviously, $Q$ varies in a system of dimension 9 on $G$, and $G$ varies
is a system of dimension 4 in $\IP^4$.  Therefore,
$$\cod(M_{8,5},M_8)=(5\times8+1)-(15+9+4)=13>5.$$ Thus the lemma holds
when $d=8$.

Assume $d=9$.  There are five possible cases: \smallbreak
 \item 1 $C$ lies in no hyperplane, and has no 7-secant line;
 \item 2 $C$ lies in some hyperplane $G$, and has no 7-secant line;
 \item 3 $C$ lies in no hyperplane, and has a 7-secant line;
 \item 4 $C$ lies in some smooth quadric surface $Q$ (and has a
7-secant line);
 \item 5 $C$ lies in some hyperplane $G$, but in no smooth quadric
surface, and has a 7-secant line.
  \smallbreak\noindent
 Here and below, we use the definition given in the middle of \p.501 of
\cite{GLP}, and say that a line is $n$-{\bf secant} to a curve if the two
schemes contain a common subscheme of length at least $n$.  Now,
 each of these five cases will be considered in turn.

In Case (1), the theorem on \p.492 of \cite{GLP} yields that
$h^1(\I_C(5))=0$.  In Case (2), Th\'eor\`eme 0.1 on \p.30 of \cite{dA}
yields the same vanishing.  Thus, in either case, $C\in M_{9,0}$.

In Case (3), the table on \p.504 of \cite{GLP} says that
$h^1(\I_{C/G}(5))=1$.  On the other hand, $C$ belongs to the subset of
$M_d$ of curves with 7-secant lines, and this subset has codimension at
least 8 by Lemma (2.4) below.  Thus the curves that fall into Case (3)
form a subset of $M_{d,1}$ of codimension at least 8 in $M_d$.

In Case (4), since $C$ is smooth, rational, and of degree 9, it is of
type (8,1) on $Q$.  So $I_{C/Q}(5)=\O_Q(-3,4)$.  So the K\"unneth
formula yields
		 $$h^1(\I_{C/Q}(5))=0\times0+2\times5=10.$$
 The linear span of $Q$ is a hyperplane $G$.  Moreover,
	$$h^1(\I_{C/Q}(5))=h^1(\I_{C/G}(5))=h^1(\I_C(5));$$
  indeed, the second equality is \Cd1a), and the first can be proved
similarly.  Now, $C$ varies in a system of dimension 17 on $Q$.  Hence,
the various $C$ in Case (4) form a subset $S_4$ of $M_{9,10}$, and
	$$\cod(S_4,M_9)=(5\times9+1)-(17+9+4)=16>10.$$

In Case (5), by Lemma (2.4) below, the various $C$ in a fixed
hyperplane $G$ form a subset of codimension at least 3 in the Hilbert
scheme of smooth rational curves of degree 9 in $G$.  Since this
Hilbert scheme has dimension $4\times9$, or 36, and since $G$ varies in
a system of dimension 4, the various $C$ form a subset of $M_9$ of
codimension at least $46-(36+4)+3$, or 9.  So it suffices to fix $C$
and to prove that $h^1(\I_C(5))<9$, or equivalently by \Cd1a), that
$h^1(\I_{C/G}(5))<9$.

Choose in $G$ a plane $H$ that meets $C$ in 9 distinct points, no
three of which are collinear.  Such an $H$ exists by \cite{ACGH, Lem.,
\p.109} or by \cite{Rath, Lem.~1.1, \p.566} (the latter is valid in
arbitrary characteristic; the curve may be singular, but must not be
strange).  Alternatively, such an $H$ exists by Bertini's theorem and
the trisecant lemma (the latter asserts that the trisecants form, at most, a
1-parameter family; it is proved for an arbitrary reduced and
irreducible space curve in characteristic 0 in \cite{L, \p.135}, and it
follows from \cite{H, Prp.~3.8, \p.311} and \cite{H, Thm.~3.9, \p.312}
for a smooth, reduced and irreducible space curve in arbitrary
characteristic).  Let $k\ge4$.  Then the 9 points impose independent
conditions on the system of curves of degree $k$ in $H$ by \cite{ACGH,
Lem., \p.115}.  Hence, in the long exact sequence,
    $$H^0(\O_H(k))\to H^0(\O_{C\cap H}(k))\to H^1(\I_{C\cap H,H}(k))
	\to H^1(\O_H(k)),$$
 the first map is surjective.  However, the last term vanishes.  Therefore,
	$$H^1(\I_{C\cap H,H} (k))=0.$$
 Consequently, the exact sequence of sheaves,
	$$0\to\I_{C/G}(k-1)\to\I_{C/G}(k)\to\I_{C\cap H,H}(k)\to0,$$
 yields
	$$h^1(\I_{C/G}(3))\ge h^1(\I_{C/G}(4))
		\ge h^1(\I_{C/G}(5))\ge\cdots.\eqlt1$$

Consider the standard exact sequence of sheaves,
	$$0\to\I_{C/G}(k)\to \O_G(k)\to\O_C(k)\to0.$$
  Since $H^1(\O_G(k))=0$, taking cohomology yields, for $k\ge0$,
  $$h^0(\I_{C/G}(k))={k+3\choose 3} -(9k+1)+h^1(\I_{C/G}(k)).\eqlt2$$
 Now, proceeding by contradiction, assume $h^1(\I_{C/G}(5))\ge9$.  Then,
	$$h^0(\I_{C/G}(3))\ge1 \and h^0(\I_{C/G}(4))\ge7.$$
 It follows, as will be proved in the next three paragraphs, that
	$$h^0(\I_{C/G}(6))\ge31.\eqlt2a$$
 Hence, $h^1(\I_{C/G}(6))\ge2$.  However, the table on \p.504 of
\cite{GLP} says that $h^1(\I_{C/G}(6))$ is equal to 1 if $C$ has an
8-secant line and to 0 if not.  Thus there is a contradiction.

Set $H(k):=H^0(\I_{C/G}(k))$, and view its elements as homogeneous
polynomials of degree $k$ in four variables that vanish on $C$.  Assume
that $H(3)$ contains a nonzero cubic $K$.  Then $K$ is irreducible
because $C$ lies in no plane, in no quadric cone, and in no smooth
quadric surface.  Now, as $L$ ranges over the linear forms, the
products $KL$ form a 4-dimensional subspace of $H(4)$.  Assume that,
modulo this subspace, $H(4)$ contains two linearly independent quartics
$Q_1$ and $Q_2$.  Then $K$ divides no nontrivial linear combination $Q$
of the $Q_i$.  So, since $K$ is irreducible, it divides no nonzero
product $QT$ where $T$ is a quadric.  Let $X$ be the set of all
products $QT$, and $Y$ the set of all products $KK'$ where $K'$ is a
cubic.  Then, therefore, $X\cap Y=\{0\}$.

Each nonzero combination $Q$ is irreducible.  Indeed, $Q$ cannot be the
product of two quadrics since $C$ lies in no quadric surface.  Suppose
$Q=K'L$ where $K'$ is a cubic.  Then $K'$ is not a multiple of $K$, but
$K'$ vanishes on $C$ because $L$ cannot.  So $C$ is contained in the
complete intersection $C'$ defined by $K$ and $K'$.  Hence $C$ is equal
to $C'$ because both have degree $9$.  However, $p_aC'=10$ by \cite{H,
Exr.~7.2(d), \p.54}, whereas $p_aC=0$.  Thus each $Q$ is
irreducible.\par
  Therefore, if $Q_1T=Q_2T'$, then $T$ and $T'$
vanish.  Hence the products $Q_1T$ and $Q_2T'$ form two 10-dimensional
linear spaces that meet only in 0 and that lie in $X$.  Therefore, $X$
is a cone of dimension at least 11.  On the other hand, $Y$ is a linear
space of dimension 20.  Since $X$ and $Y$ meet in a single point, any
smooth variety containing both of them has dimension at least 31; in
particular, their linear span does.  Thus \Cd2a) holds.  Thus
Lemma~\Cs3) is proved.
	\medbreak
To complete the proof of Proposition~\Cs2), we have to prove the
following lemma.

\lem4 Fix $d\ge4$ and $s\ge3$.  Fix $b$ with $d>b\ge3$.  In the Hilbert
scheme of rational, smooth and irreducible curves of degree $d$ in
$\IP^s$, form the subset of curves with a $b$-secant line.  Then this
subset has codimension at least $(s-1)(b-2)-b$.
 \endproclaim

Indeed, consider the space $M$ of parametrized embeddings of $\IP^1$
into $\IP^s$ with degree $d$.  Let $\Gamma$ in $\IP^1\times\IP^s\times
M$ be the family of graphs.  Then the projection $\Gamma\to\IP^1\times
M$ is an isomorphism; so $\Gamma$ is flat over $M$.  On the other hand,
the projection $\Gamma\to\IP^s\times M$ is an embedding, because its
fibers over $M$ are embeddings.  Hence, this embedding of $\Gamma$
defines a map of $M$ into the Hilbert scheme in question, and the map
is obviously surjective.  Since codimension can only decrease under
pullback along a surjection, we may replace the Hilbert scheme by $M$.

Say that the coordinate tuples on $\IP^1$ are $(t,u)$, on $\IP^s$ are
$(X_0,\dots,X_s)$, and on $M$ are $(a_{i,j})$. Then the maps are
defined by the  $s+1$ equations,
	$$X_i=a_{i,0}t^d+a_{i,1}t^{d-1}u+\cdots+a_{i,d}u^d.$$
 Fix a divisor $D$ of degree $b$ on $\IP^1$.  Say $D=P_1+\cdots+P_b$
and $P_k=(t_k,u_k)$.  Then, for each $k$, the above equations become
	$$X_i=L_{i,k}(a_{i,0},a_{i,1},\dots,a_{i,d})$$
 where the $L_{i,k}$ are linear forms.  If $P_k$ is repeated, then
replace the $n$th occurrence of $L_{i,k}$ by its $(n-1)$st derivative.

 For each $i$, the forms $L_{i,k}$ for $1\le k\le b$ are independent
because $D$ imposes $b$ conditions on the forms of degree $d$ in
$(t,u)$.  Since the variables $a_{i,j}$ change with $i$, the forms
$L_{i,k}$ for all $i,k$ are independent.  Now, consider the condition
that $D$ is carried to a length-$b$ subscheme of a line in $\IP^s$; it
becomes the condition $\rank [L_{i,k}]\le2$.  Since the $L_{i,k}$ are
independent, they may be viewed as defining a linear map from the
affine space of $s+1$ by $d+1$ matrices onto the affine space of $s+1$
by $b$ matrices.  Hence the condition defines a subset of codimension
precisely $(s-1)(b-2)$ in $M$.  As $D$ varies, the corresponding
subsets sweep out a subset of codimension $(s-1)(b-2)-c$  where $c\le b$.
Thus the lemma is proved.

\cor5 Let $F$ be a general quintic threefold in $\IP^4$, and let $C$ be
a rational, smooth and irreducible curve of degree $d\le9$ on $F$.
 \part 1 Then $C$ is embedded in $\IP^4$ with maximal rank.
 \part 2 Form the restriction to $C$ of the twisted sheaf of
differentials of $\IP^4$.  Then this locally free sheaf has\/ {\rm generic}
splitting type; namely, if $d=4n+m$ where $0\le m<4$, then
	$$\Omega^1_{\IP^4}(1)|C =\O_C(-n-1)^m\oplus\O_C(-n)^{4-m}.$$
 \part 3 If $d\le4$, then $C$ is a rational normal curve of degree $d$,
and if $4\le d\le9$, then $C$ spans $\IP^4$.
 \part 4 If $d=1$, then $C$ is $1$-regular; if $2\le d\le4$, then $C$
is $2$-regular; if $5\le d\le7$, then $C$ is $3$-regular; and if $8\le
d\le9$, then $C$ is $4$-regular.
 \endproclaim

Indeed, first let's make a general observation.  Given an arbitrary
proper closed subset $S$ of $M_d$, its preimage in $I_d$ is proper and
so has dimension at most 124 thanks to Proposition~(2.2).  Hence this
preimage does not project onto the $\IP^{125}$ of quintics.  So, since
$F$ is general, $C\notin S$.

To prove (2), in $M_d$ form the subset $S$ of curves such that the
sheaf in question does not have generic splitting type.  Then $S$ is a
proper closed subset by a theorem of Verdier's \cite{Verdier, Thm.,
p.~139} (see also \cite{Ramella, Thm.~1, p.~181}).  So $C\notin S$.
Thus (2) holds.

To prove (3), in $M_d$ form the subset $S$ of curves not spanning
$\IP^4$.  Clearly, $\dim S\le4d+4$.  If $d\ge4$, then $4d+4<5d+1$, and
so $S$ is proper.  Thus the second assertion of (3) holds.  The first
assertion follows if $d=4$, and is easy to check directly if $d\le3$
(see \cite{H, Exr.~3.4, \p.315}).

To prove (1) and (4), recall that, by definition, $C$ is embedded in
$\IP^4$ with maximal rank if, for all $k\ge1$,  the natural map,
	$$H^0(\IP^4,\O_{\IP^4}(k))\to H^0(C,\O_C(k)),$$
 is injective or surjective (or both).  It follows from the standard long
exact sequence of cohomology that (1) implies (4).

 First, suppose $d\le4$.  Then $C$ is a rational normal curve of degree $d$ by
(3).  So (1) holds because $C$ is essentially the $d$-uple image of
$\IP^1$; alternatively, (1) holds by the theorem on \p.492 of
\cite{GLP}.  Then (4) follows; alternatively, (4) holds by Remark (1)
on \p.497 of \cite{GLP}.

Finally, suppose $d\ge4$.  In $M_d$ form the subset $S$ of curves that
either don't span $\IP^4$ or aren't of maximal rank.  Then $S$ is
proper by the proof of (3) and by the maximal-rank theorem \cite{BE,
Thm.~1, \p.215}.  So $C\notin S$.  Thus (1) holds, and (4) follows.
(Curiously, (4) also follows from (2) and from Proposition~1.2 on
\p.494 of \cite{GLP}, except when $d=6,7,9$.  When $d=6,7$,
only 4-regularity comes out this way, and when $d=9$, only 5-regularity
does.)

\newsect Singular curves

In this section, we'll treat singular curves by proving the following
theorem, which complements Theorem~(2.1).  In the next section, we'll
treat reducible curves.

\thm1 On a general quintic hypersurface in $\IP^4$, there are no
rational and singular, reduced and irreducible curves of degree at most
$9$, other than the six-nodal plane quintics.
 \endproclaim

 Indeed, let $M^{r,g}_d$ be the space of all parametrized mappings of
$\IP^1$ into $\IP^4$ of degree $d$ with these (locally closed)
properties: the mappings are birational onto their images $C$, and each
such $C$ spans an $r$-plane and has arithmetic genus $g$.  Note that
each $C$ is the image of every mapping in a four parameter family ---
namely, the orbit under the parametrized automorphisms of $\IP^1$ (in
fact, the orbit includes every mapping with image $C$ because of the
hypothesis of birationality).  Let $I^{r,g}_d$ denote the incidence
locus in $M^{r,g}_d\times\IP^{125}$; it's points represent the pairs
consisting of a mapping and a quintic such that the image of the former
lies in the latter.  Except for $(r,d,g)=(2,5,6)$, the theorem asserts,
in other words, that for $r=2,3,4$ the projection
$I^{r,g}_d\to\IP^{125}$ is not surjective if $d\le9$ and $g\ge1$ (if
$g=0$, then $C$ is smooth).  Since the projection's fibers all have
dimension at least 4, to prove the that it's not surjective, we need
only establish the bound $\dim I^{r,g}_d<129$.

Consider the other projection, $I^{r,g}_d\to M^{r,g}_d$.  Its fiber
over a mapping with image $C$ is a projective space of dimension
$h^0(\I_C(5))-1$.  Let $M^{r,g}_{d,i}$ denote the (locally closed)
subset of $M^{r,g}_d$ where $h^1(\I_C(5))=i$, and let $I^{r,g}_{d,i}$
denote its inverse image in $I^{r,g}_d$.  Then, mutatis mutandis, the
proof of Equation~\Co2-1b) yields
	$$\dim I^{r,g}_{d,i} \le \dim M^{r,g}_{d,i} + (126
		-(5d+1-g)+i) -1.$$
 Therefore, to prove the theorem, it suffices to establish the bound,
	$$\dim M^{r,g}_{d,i} < 5d+5-g-i,\eqlt1$$
 unless $r=2$ and $d=5$, when a little more analysis is required.  The
proof proceeds via a case-by-case analysis.  The case where $d=9$ and
$r=3$ is special, and will be handled after the proof of Lemma~(3.4) is
complete.  So, for now, assume either that $d\le8$ or that $d=9$ and
$r\neq3$.

Since $C$ is not a rational normal curve, necessarily $d>r$.  Now,
$h^1(\I_C(5))=0$.  Indeed, if $C$ lies in a hyperplane $G$, then
$h^1(\I_C(5))$ is equal to $h^1(\I_{C/G}(5))$ by \Co2-1a), whose proof
doesn't require $C$ to be smooth.  A similar argument shows further
that, if $C$ lies in a plane $H$ contained in $G$, then
$h^1(\I_{C/G}(5)))$ is equal to $h^1(\I_{C/H}(5)))$; moreover, the
latter group always vanishes by Serre's theorem.  Suppose $C$ doesn't
lie in a plane; that is, $r=3,4$.  Then the desired vanishing holds by
the table on \p.504 of \cite{GLP}.  Thus $M^{r,g}_{d,0}=M^{r,g}_d.$

Suppose $r=2$.  Given a plane $H$, consider the subset of $M^{2,g}_d$
of mappings whose image $C$ spans $H$.  Obviously, this subset has
dimension at most $3(d+1)$.  Now, $H$ varies in a system of dimension
$3\times(5-3)$, or 6.  So
 $$\dim M^{2,g}_d\le 3d+9.$$
 On the other hand, $g=(d-1)(d-2)/2$.  Hence the bound~\Cd1) holds for
$d=3,4$.  However, for $d=5$, the two sides of \Cd1) are equal; in
other words, $I^2_{5,6}$ has dimension 129.  Now, the image $C$ of a
general mapping in $M^{2,g}_d$ has $g$ nodes and no other
singularities.  Therefore, a general quintic threefold $F$ can contain no
rational plane quintic $C$ other than one with six nodes (in fact,
Vainsencher \cite{V, \S6.2} proved there are 17,601,000 six-nodal $C$
in $F$).  Finally, suppose $d\ge6$.  Then $F$ can contain no image $C$.
Otherwise, Bezout's theorem implies that $F$ contains the span $H$ of
$C$, so infinitely many lines.

Suppose $r=3$.  Then $4\le d\le8$.  Given a hyperplane $G$, consider
the subset of $M^{3,g}_d$ of mappings whose image $C$ spans $G$.
Obviously, this subset has dimension at most $4(d+1)-1$, or $4d+3$, as
$g\ge1$.  Now, $G$ varies in a system of dimension 4.  So
 $$\dim M^{3,g}_d\le 4d+7.$$
 Hence the bound \Cd1) holds if $g<d-2$.  However, as $d\le8$, the
Castelnuovo--Halphen bounds yield $g<d-2$, except in these
seven cases:
   $$(d,g)=(6,4),\ (7,6),\ (8,8),\ (8,9), \ (7,5),\ (8,6),\ (8,7).$$
 In the first four cases, each $C$ necessarily lies on a quadric
surface.  (As is well known, these bounds and exceptions were first
proved rigorously by Gruson and Peskine \cite{GP, Thm.~3.1, \p.49},
although they assert the results only for smooth curves.  The
smoothness enters via an appeal on \p.51 to a preliminary version of
Laudal's \cite{L, Cor.~2, \p.147}; however, the final version applies
to an arbitrary reduced and irreducible curve.) Hence, Lemma~\Cs2)
below yields
	$$\dim M^{3,g}_d\le (2d+12)+4=2d+16,$$
 because the hyperplanes in $\IP^4$ form a 4-parameter family.  Thus
\Cd1) holds in those four cases.  The remaining three cases are covered
by Parts~(c) and (d) of Lemma~\Cs3) below because, if $g\ge3$, then $C$
has either at least three distinct singularities or at least one
singularity with $\delta$-invariant at least 2.

Suppose $r=4$.  Then $d\ge5$. If $d=7,8,9$, then the
Castelnuovo--Halphen bounds yield $g\le7$.  (See \cite{Cil, (1.1),
\p.27}; on \pp.24--5, Ciliberto says that the extension of the bounds
to $\IP^r$ for $r\ge4$ is due to Fano and Harris.)  Hence the bound
\Cd1) holds by Lemma~\Cs4)(a).  So suppose $d\le6$.  Then the
Castelnuovo--Halphen bounds yield $g\le2$.  If $g=1$, then \Cd1) holds
by Lemma~\Cs3)(a).  If $g=2$, then $C$ has either two distinct
singularities or one singularity of $\delta$-invariant 2; hence, \Cd1)
holds in the first case by Lemma~\Cs3)(b), and it holds in the second
by Lemma~\Cs3)(d).

Thus, if $d\le8$ or if $d=9$ and $r\neq3$, then Theorem~(3.1) holds
given Lemmas~(3.2), (3.3) and (3.4).

\lem2 Fix an integer $d$, and form the space $J$ of all pairs
$(\phi,Q)$, where $\phi$ is a parametrized generic embedding of\/
$\IP^1$ into $\IP^3$ and where $Q$ is a reduced and irreducible quadric
surface such that the image $C$ of $\phi$ is of degree $d$ and lies on
$Q$.  Then $J$ is equidimensional  of dimension $2d+12$.
 \endproclaim

Indeed, the surfaces $Q$ are parametrized by an open subset $U$ of
$\IP^9$, and those that are smooth, by a smaller open subset $V$.  The
natural action of $GL(4)$ on $\IP^3$ induces an action on $U$ and a
compatible action on $J$.  Moreover, $V$ is an orbit.  Hence, every
component of $J|V$ projects onto $V$.  Therefore, to prove the lemma,
it suffices to prove (1) that no component of $J$ projects into $U-V$
and (2) that every fiber of $J|V\to V$ is equidimensional of dimension
$2d+3$.  On the other hand, just to bound the dimension of $J$, which
is all that's needed in the proof of Theorem.~\Cs1), it suffices to
prove (2) and (3), where (3) asserts that every fiber of $J$ over a
point of $U-V$ is of dimension at most $2d+3$.

To prove (1), it clearly suffices to prove this: let $(\phi_0,Q_0)$ be
a pair where $Q_0$ is a cone, and let $T$ be the spectrum of the power
series ring in one variable $t$; then there exists a map $T\to J$ which
carries the closed point of $T$ to $(\phi_0,Q_0)$ and which carries the
generic point to a pair $(\phi',Q')$ where $Q'$ is smooth.  Say that
the coordinates on $\IP^1$ are $(u,v)$ and that those on $\IP^3$ are
$(X_0,\dots,X_3)$.  Say that
	$$\phi_0(u,v)=(\phi_{00}(u,v),\dots,\phi_{03}(u,v)),$$
 where the $\phi_{0j}$ are forms of degree $d$, and that
	$$Q_0:q_0(X_0,\dots,X_3)=0$$
 where $q_0$ is a form of degree 2.  To construct $T\to J$, it is
necessary to construct, for each $i\ge1$, similar terms,
	$$\phi_i(u,v)=(\phi_{i0}(u,v),\dots,\phi_{i3}(u,v))$$
 and $q_i(X_0,\dots,X_3)$, such that, for each $n\ge1$,
	$$(q_0+tq_1+\cdots+t^nq_n)(\phi_0+t\phi_1+\cdots+t^n\phi_n)
		\equiv0\mod(t^{n+1}).\eqlt2$$
 Of course, $q_0(\phi_0)=0$ because $(\phi_0,Q_0)\in J$.

The proof proceeds by induction on $n$.  Suppose that $\phi_i$ and
$q_i$ have been constructed for $i<n$ and that, moreover, the
$\phi_{ij}$ are linear combinations of the $\phi_{0j}$.  In the
expansion of the left hand side of \Cd2), the coefficient of $t^i$
vanishes for $i<n$ by induction.  The coefficient of $t^n$ has the form,
    $$((\grad q)(\phi_0))\cdot\phi_n+a_n(\phi_0)+q_n(\phi_0),\eqlt3$$
 where the $\grad q$ is the gradient and where $a_n$ is a certain
quadratic form.  Note that $a_1$ vanishes and that each $a_n$ is
independent of $\phi_n$ and $q_n$.  Now, $GL(4)$ acts transitively on
the space $U-V$ of cones.  So we may assume that
	$$q_0=X_0^2/2+X_1X_2.$$
 Then the first term of \Cd3) becomes
	$$\phi_{00}\phi_{n0}+\phi_{02}\phi_{n1}+\phi_{01}\phi_{n2}.$$
 Hence, if $n=1$, then \Cd3) vanishes with these choices,
  $$q_1=-(X_0^2+X_1^2+X_2^2+X_0X_3),\ \phi_{10}=\phi_{00}+\phi_{03},\
\phi_{11}=\phi_{02},\ \phi_{12}=\phi_{01},$$
 and with $\phi_{13}$ arbitrary.  If $n\ge2$, then collect the terms in
$a_n(\phi_0)$ so that it acquires the form,
 $$\phi_{00}\gamma_0+\phi_{01}\gamma_1+\phi_{02}\gamma_2+c_n\phi_{03}^2,$$
 where the $\gamma_k$ are suitable linear combinations of the
$\phi_{0j}$ and where $c_n$ is a suitable constant.  Then \Cd3) vanishes
with these choices,
  $$q_n=-c_nX_3^2,\ \phi_{n0}=-\gamma_0,\ \phi_{n1}=-\gamma_2,\
\phi_{n2}=-\gamma_1,$$
 and with $\phi_{n3}$ arbitrary.  The construction is now complete.

Set $\phi':=\sum_{n\ge0}t^n\phi_n$ and $Q':\sum_{n\ge0}t^nq_n=0$.  Then
$\phi'$ is a generic embedding because $\phi_0$ is; indeed, $\phi'$ is
finite because $\phi_0$ is, and the comorphism of $\phi'$ is
generically surjective because that of $\phi_0$ is.  By construction,
the image of $\phi'$ is of degree $d$ and lies on $Q'$.  Finally, $Q'$
is smooth because, if $c:=2(\sum_{n\ge2}t^nc_n)$, then
$\grad\sum_{n\ge0}t^nq_n$ is equal to
	$$\displaylines{
 \qquad(X_0,\ X_2,\ X_1,\ 0)-t(2X_0+X_3,\ 2X_1,\ 2X_2,\ X_0)
	-c(0,\ 0,\ 0,\ X_3)\hfill\cr
 \hfill=((1-2t)X_0-tX_3,\ X_2-2tX_1,\ X_1-2tX_2,\ -tX_0-cX_3),\qquad\cr}$$
 and the right hand side is clearly nonzero for any power series
$X_i(t)$ and $c(t)$, provided at least one $X_i$ is nonzero and $c(t)$
is not $t^2/(1-2t)$; the latter restriction can be achieved by an
appropriate choice of $\phi_{13}$.  Thus (1) holds.

To prove (2), set $Q:=\IP^1\times\IP^1$.  Then a map $\phi\:\IP^1\to Q$
is given by a pair of maps $\phi_i\:\IP^1\to\IP^1$ for $i=1,2$.  Set
$d:=\deg\phi$ and $d_i:=\deg\phi_i$.  Then $d=d_1+d_2$.  Say that the
coordinates on $\IP^1$ are $(u,v)$ and that those on $\IP^3$ are
$(X_0,\dots,X_3)$.  Then a parametrization,
    $$(\alpha_i(t,u),\beta_i(t,u)),$$
 of each $\phi_i$ yields a parametrization of $\phi$, namely,
	$$\bigl(\alpha_1(t,u)\beta_1(t,u),\alpha_1(t,u)\beta_2(t,u),
\alpha_2(t,u)\beta_1(t,u),\alpha_2(t,u)\beta_2(t,u)\bigr).$$
 If the parametrization of $\phi_i$ is multiplied by a scalar $v_i$,
then the parametrization of $\phi$ is multiplied by $v_1v_2$.  Form
the space of all pairs of parametrized maps $\IP^1\to\IP^1$ of bidegree
$(d_1,d_2)$; obviously, this space has dimension,
	$$2(d_1+1)+2(d_2+1)=2d+4.$$
 Now, it is an open condition on a map $\phi\:\IP^1\to Q$ to be a
generic embedding.  Therefore, each fiber of $J|V\to V$ breaks up into
the disjoint union of subspaces, each of dimension $2d+3$.  Thus (2)
holds.  Finally, (3) was proved by Guest \cite{Gu, Prop.~2.1, \p.60} via
a similar, but a little more complicated, elementary description of the
parametrized maps from $\IP^1$ to the cone with equation $X_2^2-X_1X_3$.
The proof of Lemma~\Cs2) is now complete.

\lem3 Fix $d\ge5$, and let $M^r_{d,n,\delta}$ be the space of
parametrized mappings of $\IP^1$ into $\IP^4$ with these (locally
closed) properties: each mapping is birational onto its image $C$, and
this $C$ spans an $r$-plane, has degree $d$, and has at least $n$
distinct singularities of which one has $\delta$-invariant at least
$\delta$.
 \part a Then $\dim(M^4_{d,1,1})\le5d+3$ and $\dim(M^3_{d,1,1})\le4d+7$.
 \part b Then $\dim(M^4_{d,2,1})\le5d+1$ and $\dim(M^3_{d,2,1})\le4d+6$.
 \part c Then $\dim(M^4_{d,3,1})\le5d-1$ and $\dim(M^3_{d,3,1})\le4d+5$.
 \part d Then $\dim(M^4_{d,1,2})\le5d$ and $\dim(M^3_{d,1,2})\le4d+5$.
 \part e Then $\dim(M^4_{d,2,2})\le5d-2$ and $\dim(M^3_{d,2,2})\le4d+4$.
 \part f Then $\dim(M^4_{d,1,3})\le5d-2$ and $\dim(M^3_{d,1,3})\le4d+4$.
 \pf
 Say that the coordinate tuples on $\IP^1$ are $(t,u)$, on $\IP^4$ are
$(X_0,\dots,X_4)$, and on $M^r_{d,n,\delta}$ are $(a_{i,j})$.  Then the
maps are defined by the five equations,
	$$X_i=a_{i,0}t^d+a_{i,1}t^{d-1}u+\cdots+a_{i,d}u^d.\eqlt3a$$
 We'll treat the case $r=4$ only, as the case  $r=3$ is similar.

To prove (a), consider the condition that two (distinct) points of
$\IP^1$ are mapped to the same point of $\IP^4$.  If these points are
$(1,0)$, $(0,1)$, and $(1,0,0,0,0)$, then the condition becomes this:
$a_{i,0}=0$ and $a_{i,d}=0$ for $i=1,2,3,4$.  Hence, the corresponding
space of mappings is of dimension $5d-3$.  Let the first two points
vary over $\IP^1$, and the third over $\IP^4$.  Then the corresponding
spaces sweep out a space of dimension at most $5d+3$.  Now, consider
the condition that a point of $\IP^1$ is mapped to a cusp on $C$.  If
these points are $(1,0)$ and $(1,0,0,0,0)$, then the condition becomes
this: $a_{i,0}=0$ and $a_{i,1}=0$ for $i=1,2,3,4$.  Hence, the
corresponding space of mappings is of dimension $5d-3$ too.  Let the
first point vary over $\IP^1$, and the second over $\IP^4$.  Then the
corresponding spaces sweep out a space of dimension at most $5d+2$.
Thus (a) holds.

To prove (b), consider the condition that four points of $\IP^1$ are
mapped 2-to-1 to two points of $\IP^4$.  If the four are $(1,0)$,
$(0,1)$, $(1,1)$ and $(1,v)$ with $v\neq0,1$ and if the first two map
to $(1,0,0,0,0)$ and the second two map to $(0,0,0,0,1)$, then the
condition becomes this: for $i=1,2,3,4$,
	$$a_{i,0}=0,\ a_{i,d}=0,\ \sum_ja_{i-1,j}=0,\
\sum_ja_{i-1,j}v^j=0.\eqlt4$$
 These sixteen linear equations are obviously independent.  Hence, the
corresponding space of mappings is of dimension $5d-11$.  Let the six
points vary.  Then the corresponding spaces sweep out a space of
dimension at most $4+4+4$ more, or at most $5d+1$ in all.

Consider the condition that a point of $\IP^1$ is mapped to a cusp and
that a pair of points is mapped to a node.  If the first point is
$(1,0)$ and the pair consists of $(1,1)$ and $(1,v)$ with $v\neq0,1$
and if the cusp lies at $(1,0,0,0,0)$ and the node at $(0,0,0,0,1)$,
then the condition is expressed by sixteen independent linear
equations; they are nearly the same as those displayed in \Cd4), only
$a_{i,d}=0$ is replaced by $a_{i,1}=0$.  Let the five points vary.  Then
the corresponding spaces of mappings sweep out a space of dimension at
most $5d$.  Similarly, the condition that two points of $\IP^1$ are
mapped to distinct cusps defines a space of mappings of dimension at
most $5d-1$.  Thus (b) holds.

In (c), there are essentially two different cases: (i) the three points
on $C$ are not collinear, say they are
	$$(1,0,0,0,0),\ (0,0,0,0,1),\ (0,1,0,0,0);$$
 and (ii) the three points are collinear, say they are
	$$(1,0,0,0,0),\ (0,0,0,0,1),\ (a,0,0,0,b)$$
 for suitable $a$, $b$.  Suppose first that all three points are nodes
of $C$.  Say the first is the image of $(1,0)$ and $(0,1)$; the second,
of $(1,1)$ and $(1,v)$; and the third of $(1,w)$ and $(1,x)$.  Then,
as above, the sixteen independent linear equations \Cd4) hold.  In Case
(i), there are these eight additional equations: for $i=0,2,3,4$,
	$$\sum_ja_{i,j}w^j=0\and \sum_ja_{i,j}x^j=0.$$
  All twenty-four equations are independent since $d\ge5$.  Hence they
define a space of mappings of dimension $5d-19$.  Let the nine points
vary.  Then the corresponding spaces sweep out a space of dimension at most
$6+3\times4$ more, or at most $5d-1$ in all.  In Case (ii), six of the
additional equations hold~--- namely, those for $i=1,2,3$~--- and all
twenty-two equations are independent.  Hence they define a space of
mappings of dimension $5d-17$.  However, the third point in $\IP^4$ has
only one degree of freedom, not four, given the other two.  Hence, as
the nine points vary, the corresponding spaces sweep out a space of
dimension at most $6+2\times4+1$ more, or at most $5d-2$ in all.  Suppose
now that one of the three points on $C$ is a cusp.  Then the argument
can be modified as in (a) or (b) to yield asserted bound.  Now, if $C$
has two cusps (or more), then the asserted bound holds by the argument
in (b).  Thus (c) holds.

To prove (d), note that a singularity on $C$ with $\delta$-invariant at
least 2 falls into one of four cases: it has at least
 \smallskip
 \item i three branches, or

 \item ii two branches, one with a cusp, or

 \item ii$'$ two branches with the same tangent line, or

 \item iii one branch with a higher cusp.
 \smallskip

Cases (i) and (ii) can be handled much as in (b).  Indeed, consider (i).  If
$(1,0)$, $(0,1)$, and $(1,1)$ map to $(1,0,0,0,0)$, then, for $i=1,2,3,4$,
	$$a_{i,0}=0,\ a_{i,d}=0,\ \sum_ja_{i,j}=0.\eqlt5$$
 These twelve linear equations are obviously independent.  Hence, the
corresponding space of mappings is of dimension $5d-7$.  Let the four
points vary.  Then the corresponding spaces sweep out a space of
dimension at most $3+4$ more, or at most $5d$ in all.

Consider (ii).  If $(1,0)$ is mapped to a cusp at $(1,0,0,0,0)$ and if
$(0,1)$ is also mapped to $(1,0,0,0,0)$, then, for $i=1,2,3,4$,
	$$a_{i,0}=0,\ a_{i,1}=0,\ a_{i,d}=0.$$
 These twelve linear equations are obviously independent.  Hence, the
corresponding space of mappings is of dimension $5d-7$.  Let the three
points vary.  Then the corresponding spaces sweep out a space of dimension
at most $2+4$ more, or at most $5d-1$ in all.

Consider (ii$'$).  If $(1,0)$ and $(0,1)$ are mapped to $(1,0,0,0,0)$,
then $a_{i,0}=0$ and $a_{i,d}=0$ for $i=1,2,3,4$.  If also the
corresponding tangent lines are equal, then there is a scalar $r$ such
that $a_{i,d-1}=ra_{i,1}$ for $i=0,\dots,4$.  Hence, the corresponding
space of mappings is of dimension $5d-7$.  Let the three points vary.
Then the corresponding spaces sweep out a space of dimension at most
$5d-1$.

Finally, consider (iii).  Say $(1,0)$ maps to a higher cusp at
$(1,0,0,0,0)$.  Then $a_{i,0}=0$ and $a_{i,1}=0$ for $i=1,2,3,4$.  In
addition, if the cusp has multiplicity at least 3, then $a_{i,2}=0$ for
$i=1,2,3,4$; if not, then no linear combination of the $X_i$ can vanish
to order 3 at $(1,0)$, and therefore there is a scalar $r$ such that
$a_{i,3}=ra_{i,2}$ for $i=1,\dots,4$.  Hence, the corresponding space of
mappings is of dimension $5d-6$.  Let the two points vary.  Then the
corresponding spaces sweep out a space of dimension at most $1+4$ more,
or at most $5d-1$ in all.  The proof of (d) is now complete.

The proofs of (e) and (f) build on the proof of (d).  The latter gives
more than asserted; namely, in Cases~(ii), (ii$'$), and (iii), the
space of mappings is irreducible of dimension at most $5d-1$.  If the
dimension is, in fact, at most $5d-2$, then we're done.  Assume not.
Then, in each case, a general curve has exactly one singularity of
$\delta$-invariant exactly 2 and no singularity of greater
$\delta$-invariant, as we'll now show.

In Case~(ii), proceed as follows.  For $r=4$, a general projection from
the point $P:=(1,0,0,0,0)$ gives a smooth curve of degree $d-3$ in
$\IP^3$; indeed, simply disregard the common factor of $tu^2$ in the
parametrization of the projected curve.  Hence the original curve is
smooth outside $P$, and at $P$ has no other branches than $t=0$ and
$u=0$.  Moreover, the curves with $\delta$-invariant at least 3 at $P$
and no other branches at $P$ clearly form a proper, closed subset.
Indeed, their parametrizations satisfy these additional conditions: if
there's a higher cusp at $u=0$, then $ka_{i,2}+la_{i,3}=0$ for
$i=1,\dots,4$ and for some fixed $k,\ l$; if there's a cusp at $t=0$,
then $a_{i,d-1}=0$ for $i=1,\dots,4$; and if there's a common tangent
at the two branches, then $ka_{i,2}+la_{i,d-1}=0$ for $i=1,\dots,4$.
For $r=3$, the argument is nearly identical.  However, the general
projection from $(1,0,0,0)$ yields a plane curve of degree $d-3$.  So
it's not smooth if $d\ge 6$, but has ordinary nodes.  Since the
parametrization of $X_0$ is arbitrary, the original curve is smooth
nevertheless.

In Cases~(ii$'$) and (iii), we argue similarly.  In these cases, the
general parametrizations give curves, which are smooth along $P$, and
have no branches passing through $P$, except for $t=0$ and $u=0$.  Hence
the special parametrizations giving a singularity of $\delta$-invariant
at least 3 at $P$ form a proper closed subset.

Thus in Parts~(e) and (f) the mappings of Cases~(ii), (ii$'$) and (iii)
form proper closed subspaces.  Hence their dimensions are at most
$5d-2$.  So, it remains to handle just the mappings of Case~(i),
those whose image $C$ has at least three branches at one point.
(Alternatively, the previous three cases can be handled similarly.)
 If, say, $(1,0)$, $(0,1)$, and $(1,1)$ map to $(1,0,0,0,0)$, then the
twelve equations \Cd5) are satisfied.

In (e), each $C$ has an additional singularity.  If this singularity
lies at $(0,0,0,0,1)$ and if it's the image of $(w,1)$ and $(x,1)$, then
the following eight additional equations are satisfied: for $i=0,1,2,3$,
	$$\sum_ja_{i,j}w^{d-1-j}=0\and \sum_ja_{i,j}x^{d-1-j}=0.$$
  All twenty equations are independent, by the theory of generalized
Vandermonde determinants, if $w$ and $x$ are general since $d\ge5$.
Hence they define a space of mappings of dimension $5d-15$.  Let the
seven points vary.  Then the corresponding spaces sweep out a space of
dimension at most $5+2\times4$ more, or at most $5d-2$ in all.  If
instead $(0,0,0,0,1)$ is a cusp, and it's the image of $(w,1)$, then
the following eight additional equations are satisfied: for $i=0,1,2,3$,
	$$\sum_ja_{i,j}w^{d-1-j}=0\and \sum_j(d-j)a_{i,j}x^{d-1-j}=0.$$
  Again, all twenty equations are independent, but this time there are
only six points to vary.  Hence the space of maps is of dimension
at most $5d-3$.

In (f), each $C$ has $\delta$-invariant at least 3 at the singular
point.  Two configurations are possible: \smallskip
 \item i four or more branches,
 \item ii three branches whose tangent lines are coplanar.
 \smallskip
 Consider (i).  If also $(v,1)$ maps to $(1,0,0,0,0)$, then the
following four additional equations are satisfied: for $i=1,2,3,4$,
	$$\sum_ja_{i,j}v^{d-j}=0.$$
 All sixteen equations are independent if $v\neq 0,1$.  Hence they
define a space of mappings of dimension $5d-11$.  Let the five points
vary.  Then the corresponding spaces sweep out a space of dimension at
most $4+4$ more, or at most $5d-3$ in all.

Finally, consider (ii).  If $(1,0)$, $(0,1)$, and $(1,1)$ map to
$(1,0,0,0,0)$, then, since their tangent lines are coplanar, the maximal
minors of the following matrix vanish:
	$$\pmatrix{a_{11}&a_{21}&a_{31}&a_{41}\cr
      a_{1,d-1}&a_{2,d-1}&a_{3,d-1}&a_{4,d-1}\cr
      c_1&c_2&c_3&c_4\cr}$$
 where $c_i:=\sum_{j=1}^{d-1}(d-j)a_{i,j}$.  These condition, together
with the twelve equations \Cd5), define a space of mappings of
dimension $5d-9$.  Let the four points vary.  Then the corresponding
spaces sweep out a space of dimension at most $3+4$ more, or at most
$5d-2$ in all.  The proof of Lemma~\Cs3) is now complete.

\lem4 (a) If $d\ge7$, then $\dim M^{4,g}_d\le5(d+1)-\min(2g,8)$.
 \part b If $d\ge9$, then $\dim M^{3,g}_d\le 4(d+1)-\min(g,5)$.
 \pf
 For $g\le3$, the assertions follow from Lemma~\Cs3).  In the remaining
cases, the proof is similar to that of this lemma, and the same
notation will be used.  For $g\ge4$, there are five cases of
configurations of minimal singularities on $C$, corresponding to the
partitions of 4: \smallskip
 \item i four singularities with $\delta$-invariant 1,
 \item ii one singularity with $\delta$-invariant 2 and two with
$\delta$-invariant 1,
 \item iii two singularities with $\delta$-invariant 2,
 \item iv one singularity with $\delta$-invariant 3 and one with
$\delta$-invariant 1,
 \item v one singularity with $\delta$-invariant 4.
 \smallskip\noindent
 Each case will now be considered in turn.

In Case (i), the proof is similar to the proofs of (a), (b) and (c) of
Lemma~\Cs3).  Indeed, first fix the locations of the four singularities
in $\IP^4$ (resp., in $\IP^3$) and the locations of their preimages in
$\IP^1$.  Doing so imposes $4\times8$ (resp., $4\times6$) independent
linear conditions on the coefficients $a_{i,j}$ in the equations~\Cd3a)
because $d\ge7$.  If the four singularities are nodes, then, when the
locations of the twelve points are varied, the corresponding linear
spaces of mappings sweep out a space of dimension at most $8+4\times4$
more (resp., $8+4\times3$ more), or dimension at most $5d-3$ (resp.,
$4d$) in all.  If $n$ of the singularities are cusps, then the
dimension is $n$ more.  Thus the assertions hold in Case~(i).

In Cases (ii) and (iii), the proofs are similar to those of (d) and (e)
of Lemma~\Cs3).  First consider a singularity with $\delta$-invariant
2.  Fix its type.  Fix its location and the locations of its preimages.
Recall that doing so imposes independent linear conditions on the
coefficients $a_{i,j}$ in the equations~\Cd3a).  If there are three
preimages, then, when the locations of the four points are varied, the
corresponding linear spaces of mappings sweep out a space of dimension
at most $5d$ (resp., $4d+1$).  If there are fewer preimages, then the
dimension is smaller.  Now, it is not hard to see similarly that,
because $d\ge7$, imposing a second singularity with $\delta$-invariant
2 decreases the dimension by at least 5 more.  Thus the assertions hold
in Case~(iii).  On the other hand, it is not hard to see that imposing
two additional singularities with $\delta$-invariant 1 decreases the
dimension by at least 4 (resp., 2).  Thus the assertions hold in
Case~(ii).

In Case (iv), repeat the argument proving (f) of Lemma~\Cs3), and note
that imposing an additional singularity with $\delta$-invariant 1
decreases the dimension by at least 2 (resp., 1).

In Case (v), the proof is similar to the proofs of (d) and (f) of
Lem\-ma~\Cs3).  There are five subcases, according to the number of
branches.  If the number is 5, then the space of mappings has dimension
at most $5(d+1)-5\times4+(5+4)$, or $5d-6$ (resp.,
$4(d+1)-5\times3+(5+3)$, or $4d-3$).  If the number is 4, then the
space of mappings has dimension at most $5(d+1)-4\times4+(4+4)$, or
$5d-3$ (resp., $4(d+1)-4\times3+(4+3)$, or $4d-1$); in fact, the
dimension is less because the four tangent lines must span a 3-space,
but the crude bound is sufficient.  If the number is 3, then, by the
proof of (d) of Lemma~\Cs3), the space of mappings with coplanar
tangent lines has dimension at most $5d-2$ (resp., $4d$); however,
since $\delta=4$, two of the lines must coincide or one branch must be
cuspidal, and so the space of mappings in question must have dimension
at most $5d-3$ (resp., $4d-1$).  Finally, if the number is 2 or 1,
then, by the proof of (d) of Lemma~\Cs3), the space of mappings such
that $\delta\ge2$ has dimension at most $5d-1$ (resp., $4d$); pursuing
the same line of reasoning further, it is easy to see that the subspace
of mappings such that $\delta\ge3$ has dimensions at most $5d-3$
(resp., $4d-1$).

The proof of (a) is now complete.  Moreover, in all cases of (b)
considered so far, the space of mappings has dimension at most $4d-1$,
except in the case of four singularities with $\delta$-invariant 1,
where the dimension is $4d$.  However, the reasoning in that case shows
that, in the case of five singularities with $\delta$-invariant 1, the
dimension is at most $4d-1$ because $d\ge9$.  The proof of (b) is now
complete.  Thus Lemma~\Cs4) is proved.

\medbreak

To complete the proof of Theorem~\Cs1), it remains to establish the
bound~\Cd1) for $r=3$ and $d=9$, namely, the bound,
	$$\dim M^{3,g}_{9,i} < 50-(g+i).\eqlt6$$
 Consider the image $C$ of a mapping in $M^{3,g}_9$.  Then
$i=h^1(\I_C(5))$, and as was noted before, the latter $h^1$ is equal to
$h^1(\I_{C/G}(5))$, where $G$ denotes the hyperplane spanned by $C$;
thus,
	$$i=h^1(\I_{C/G}(5)).$$

The proof of (2-1) yields, in the present case, that, for $k\ge0$,
	$$h^0(\I_{C/G}(k))={k+3\choose 3}
		-\bigl(9k+1-g+h^1(\O_C(k))\bigr)+h^1(\I_{C/G}(k)).$$
 However, $h^1(\O_C(k)))=0$ if $9k>2g-2$, because then
$\deg(\omega_C(-k))< 0$ where $\omega_C$ is the dualizing sheaf, which
is torsion free of rank 1.  Since $d=9$, the Castelnuovo--Halphen
bounds yield $g\le12$.  Hence, in particular,
	$$g+h^1(\I_{C/G}(3))=h^0(\I_{C/G}(3))+8.$$
 Now, the proof of \Co2-1) does not require $C$ to be smooth.  Therefore,
	$$g+h^1(\I_{C/G}(5))\le h^0(\I_{C/G}(3))+8.\eqlt7$$

Suppose that $C$ lies in a quadric surface $Q$.  Then $Q$ is
irreducible as $C$ lies in no plane.  If $C$ lies in a second surface
$K$ in $G$, then $C$ lies in $Q\cap K$.  So, if $\deg K\le4$, then
$K\supset Q$ by Bezout's theorem since $d=9$ and $Q$ is irreducible.
Hence, if $\deg K=3$, then $K=Q+H$ where $H$ is a hyperplane.
Therefore, $h^0(\I_{C/G}(3))=4$.  Hence, \Cd7) yields
$g+h^1(\I_{C/G}(5)) \le12$.  So the right hand side of \Cd6) is at least
38.  On the other hand, Lemma~\Cs2) implies that the various $C$ in
question correspond to a subspace $J$ of $M^{3,g}_9$ of dimension 30.
Consequently, we may assume that $C$  lies in no quadric.

Suppose that $g+h^1(\I_{C/G}(5))\ge10$.  Then \Cd7) implies that $C$
lies on two different cubics, $K$ and $G$ say.  They are irreducible
because $C$ lies in no quadric.  Hence their complete intersection is a
curve of degree 9 containing $C$, so equal to $C$.  Therefore, $g=10$
by \cite{H, Exr.~7.2(d), \p.54} for example.  Moreover,
$h^1(\I_{C/G}(5))$ vanishes; indeed, in the exact sequence of twisted
ideals,
	 $$0\to\I_{K/G}(5)\to\I_{C/G}(5)\to\I_{C/K}(5)\to0,$$
 the first term is equal to $\O_G(2)$ and the third term is equal to
$\O_K(2)$, and so both terms have vanishing $h^1$.  Thus $C$ is among
the curves treated in the next paragraph.

Suppose that $g\le 10$ and $h^1(\I_{C/G}(5))=0$.  Then the right hand
side of \Cd6) is equal to $50-g$, so is at least 40.  On the other
hand, Lemma~\Cs4)(b) implies that $M^{3,g}_9$ has dimension at most 39.

Finally, suppose that $g< 10$ and $h^1(\I_{C/G}(5))\neq0$.  Since
$d=9$, the latter condition implies that $C$ has a 7-secant line by
virtue of Theorem~0.1 on \p.30 of \cite{dA} (although the statement in
print does not make it completely clear, the theorem does indeed apply
to any singular, reduced and irreducible curve, according to a personal
communication from its author on 950705).  Since $g< 10$, by the second
paragraph above, $g+i<10$; hence, the right hand side of \Cd6) is at
least 41.  Hence \Cd6) holds by the following lemma, and once it's
proved, the proof of Theorem~\Cs1) will be complete.

\lem5 The mappings $\phi\:\IP^1\to\IP^3$ of degree $d$ at least $8$
whose images $C$ are singular and have a $7$-secant line form a space
of dimension at most $4d$.
 \pf
 The proof is similar to the last part of that of Lemma~(2.4).  Fix a
divisor $D$ of degree 7 on $\IP^1$, and fix a line $L$ in $\IP^3$.
Consider the mappings $\phi$ such that $\phi^{-1}L\supset D$.  It is
not hard to see that the $\phi$ form a linear space of dimension
$4d-10$.  The corresponding images $C$ need not have $L$ as a 7-secant
because the scheme $C\cap L$ might have length less than 7.  However,
every $C$ with $L$ as a 7-secant arises in this way.  Now, let $D$ and
$L$ vary.  Then the corresponding linear spaces sweep out an irreducible
space of dimension at most $7+4$ more, or $4d+1$ in all.  Since a
general point represents an embedding (see \cite{CJ, Prop.~4.4, \p.80}
for instance), the points representing mappings with singular images
form a subset of smaller dimension.  Thus the assertion holds.

\rmk6
 Even though a general quintic threefold does contain many six-nodal
quintic plane curves, it does not contain any rational, reduced and
irreducible, decic curve that lies on a quadric surface (possibly
singular).  Indeed, such a curve $C$ would be the complete intersection
of the quadric and the quintic, since a general quintic contains no
quadric surface.  It follows that $C$ has arithmetic genus 16 and that
$h^1(\I_C(5))$ vanishes.  Hence, by the reasoning that led to the bound
\Cd1), it suffices to prove that the various possible $C$ form a space
of dimension at most $5\times10+5-16-0$, or 39.  However, by
Lemma~(3.2), this space has dimension at most $2\times10+12+4$, or 36,
and the proof is complete.

\newsect        Reducible curves

In this section we will prove the following theorem, which completes
 our treatment of Clemens' conjecture in degree at most 9.  The theorem
says, in other words, that on a general quintic threefold there is no
pair of intersecting rational, reduced and irreducible curves whose
degrees total at most 9.

 \thm1 On a general quintic threefold in $\IP^4$, there is no
connected, reducible and reduced curve of degree at most $9$ whose
components are rational.
 \endproclaim

Indeed, suppose there is such a curve $C$.  Obviously, we may assume
$C$ has two components.  Consider one of them.  By Theorem~(3.1),
either it's a six-nodal plane quintic or it's smooth.  If it's smooth,
then, by Corollary~(2.5)(3), either it's a rational normal curve or it
spans $\IP^4$, and by Corollary~(2.5)(2), the restricted twisted sheaf
of differentials of $\IP^4$ has generic splitting type.  We are now
going to prove that, in fact, there can be no such $C$.

Let $M_a'$ be the open subscheme of the Hilbert scheme of $\IP^4$
parametrizing the smooth irreducible curves of degree $a$ that are
rational normal curves if $a\le4$ and that span $\IP^4$ if $a\ge4$; in
addition, assume that the restriction of $\Omega^1_{\IP^4}(1)$ has
generic splitting type.  Let $N_5$ be the scheme parametrizing the
six-nodal plane quintics in a variable plane in $\IP^4$.  Let
$R_{a,b,n}$ (resp., $S_{a,n}$) be the locally closed subset of $M_a'\x
M_b'$ (resp., of $M_a'\x N_5$) of pairs $(A,B)$ such that $A\cap B$ is
exactly of length $n$.  Let $I_{a,b,n}$ (resp., $J_{a,n}$) be the
locally closed subset of $R_{a,b,n}\x \IP^{125}$ (resp., of
$S_{a,n}\x\IP^{125}$) of triples $(A,B,F)$ such that $A\subset F$ and
$B\subset F$.  The $F$ that contain a plane form a proper closed subset
of $\IP^{125}$.  Form its complement, and the latter's preimage in
$J_{a,n}$.  Then replace $S_{a,n}$ by the image of that preimage, and
replace $J_{a,n}$ by the preimage of the new $S_{a,n}$.  Then, given
any pair $(A,B)$ in $S_{a,n}$, there is an $F$ that contains both $A$
and $B$, but not any plane.  It suffices now to prove that $I_{a,b,n}$
(resp., $J_{a,n}$) has dimension at most 124 for $a+b\le9$ (resp.,
$a\le4$) and $n\ge1$.

The fiber of $I_{a,b,n}$ (resp., $J_{a,n}$) over a pair $(A,B)$ is a
projective space of dimension $h^0(\I_C(5))-1$, where $C$ is the
reduced curve $A\cup B$ and $\I_C$ is its ideal in $\IP^4$.  Hence,
mutatis mutandis, the proof of Equation~\Co2-1b) yields
	$$\displaylines{\dim I_{a,b,n} \le \dim R_{a,b,n} +125
		-\min_C\{h^0(\O_C(5))-h^1(\I_C(5))\}\cr
 \bigl(\hbox{resp., }\dim J_{a,n} \le \dim S_{a,n} + 125
		-\min_C\{h^0(\O_C(5))-h^1(\I_C(5))\}\bigr).\cr}$$
 To bound the term $h^0(\O_C(5))$, consider the exact sequence,
	$$0\to\I_{B/C}\to\O_C\to\O_B\to0$$
 where $\I_{B/C}$ is the ideal of $B$ in $C$.  Now, $\I_{B/C}$ is equal
to the ideal $\I_{D/A}$ where $D:=A\cap B$.  Hence,
	$$h^0(\O_C(5))\ge\chi(\O_C(5))=\chi(\I_{D/A}(5))+\chi(\O_B(5)).$$
 Now, $A$ and $B$ are smooth, rational and of degrees $a$ and $b$
(resp., $B$ is a six-nodal plane quintic), and $D$ is of length $n$.
Hence
	$$h^0(\O_C(5))\ge 5(a+b)+2-n
	\hbox{ (resp., }h^0(\O_C(5))\ge 5a+21-n).$$
 Theorem~\Cs1) is now a consequence of the following two lemmas.

\lem2 With the notation as above, for $a+b\le9$ (resp., $a\le4$) and
$n\ge1$,
	$$\dim R_{a,b,n}\le 5(a+b)+1-n \hbox{ \rm\leftp resp., }
	\dim S_{a,n}\le 5a+20-n).$$
 \endproclaim

\lem3 With the notation as above, for $a+b\le9$ (resp., $a\le4$) and
$n\ge1$,
	$$h^1(\I_C(5))=0.$$
 \endproclaim

To prove Lemma~\Cs2), first consider a pair $(A,B)$ in $R_{a,b,n}$.  Fix
$B$ momentarily, and let $A$ vary in the fiber $\Phi_B$ of $R_{a,b,n}$
over $B$.  We may assume that $a\le b$, whence $a\le4$ and $b\le9-a$.
We'll consider each of the four values of $a$ in turn.  First, suppose
$a=1$.  Then $A$ is an $n$-secant line of $B$.  If $n=1$, then
$\dim\Phi_B=4$ since the lines through a fixed point of $B$ form a
$\IP^3$.  Now, $M_b'$ has dimension $5b+1$.  So $R_{1,b,1}$ has
dimension $5b+5$, the asserted bound.  If $n=2$, then $\dim\Phi_B=2$
since a bisecant is determined by its ``end-points.''  So $R_{1,b,2}$
has dimension $5b+3$, one less than the asserted bound.  Suppose
$n\ge3$.  Then $b\ge4$ since neither a conic nor a twisted cubic has
even a trisecant.  So $B$ spans $\IP^4$.  Hence $n\le b-2$; indeed, if
$H$ is the hyperplane spanned by an $n$-secant line $A$ and two general
points of $B$, then $H\cap B$ has length $b$.  Therefore, Lemma~(2.4)
implies that the various $B$ with an $n$-secant line $A$ form a space of
dimension at most $5b-2n+7$.  Now, $\dim\Phi_B\le2$ since $n\ge2$.  So
$R_{1,b,n}$ has dimension at most $5b-2n+9$, which is at most the
asserted bound of $5b-n+6$ since $n\ge3$.  Thus Lemma~\Cs2) holds when
$a=1$.

Suppose $a=2$.  Then $A$ is an $n$-secant conic of $B$, and we argue
much as above.  If $n=1$, then $\dim\Phi_B=9$ since there is a $\IP^4$
of conics that pass through a fixed point $P$ of $B$ and that lie in a
fixed plane and since the planes through $P$ form a 4-dimensional
space.  So $R_{2,b,1}$ has dimension $5b+10$, the asserted bound.  If
$n=2$, then $\dim\Phi_B=7$ since there is a $\IP^3$ of conics that
contain a fixed length-2 subscheme $W$ of $B$ and that lie a fixed
plane and since the planes through $W$ form a $\IP^2$.  So $R_{2,b,2}$
has dimension $5b+8$, one less than the asserted bound.  Suppose $n=3$.
Then $\dim\Phi_B\le5$ because the length-3 subschemes $W$ of $B$ form a
threefold, and the conics through $W$ form a surface when $W$
determines a plane.  So $R_{2,b,3}$ has dimension at most $5b+6$, two
less than the asserted bound for $n=3$.
 Since the same reasoning applies to both $R_{2,b,4}$ and $R_{2,b,5}$,
the latter have dimension at most $5b+6$, the asserted bound for $n=5$.
Suppose $n=5$ and $b\ge4$.  Then $\dim\Phi_B\le4$.  Indeed, there is at
most one conic through a given length-5 scheme $W$, and if there is
one, then $W$ must be planar.  Now, the length-5 subschemes $W$ of $B$
are parametrized by its fifth symmetric product, which is irreducible,
and a general such $W$ is not planar since $B$ spans $\IP^4$ as
$b\ge4$.  Thus $\dim\Phi_B\le4$.  Hence $\dim R_{2,b,5}\le5b+5$.  Since
the same reasoning applies to $R_{2,b,6}$, the latter has dimension at most
$5b+5$, the asserted bound for $n=6$ if $b\ge4$.

It now suffices to show that, if $n\ge6$, then $n=6$ and $b=7$.  Recall
that $a=2\le b$.  Denote the plane of $A$ by $G$.  If $b=2$, then $B$
is a conic, say with plane $H$.  If $H=G$, then $n=4$.  If not, set
$L:=H\cap G$.  Then $A\cap B$ is a subscheme of $A\cap L$, and the
latter has length at most 2; so $n\le2$.  If $b=3$, then $B$ is a
twisted cubic, say spanning the hyperplane $H$.  If $H\supset G$, then
$G\cap B$ has length 3, but it contains $A\cap B$. If $H\not\supset G$,
then $H\cap A$ has length 2, but it contains $A\cap B$; so $n\le2$.
Finally, if $b\ge4$, then $B$ spans the ambient $\IP^4$.  Let $H$ be
the hyperplane spanned $G$ and a general point of $B$.  Then $H\cap B$
has length $b$.  Hence $n\le b-1$, but $b\le 9-a$ and $a=2$.  Thus
$n\le6$, and if $n=6$, then $b=7$, as required.  Thus Lemma~\Cs2) holds
when $a=2$.

Suppose $a=3$.  Then $A$ is a twisted cubic, $n$-secant to $B$, and
$3\le b\le 6$.  Let $G$ and $H$ be the linear spans of $A$ and $B$.
Suppose $b=3$.  Then $B$ is a twisted cubic too.  So both $G$ and $H$
are hyperplanes.  If they are distinct, then $n\le3$, because $G\cap H$
is plane, so meets $A$ is a scheme of length $3$, and it contains
$A\cap B$.  If $G=H$, then $n\le6$, because $A$ lies in a quadric
surface $Q$ in $H$ not containing $B$, and $Q\cap B$ has length 6.
Suppose $b\ge4$.  Then $B$ spans $\IP^4$.  So $G$ meets $B$ in a
subscheme of length $b$, and it contains $A\cap B$.  Hence $n\le
b\le6$.  Now, the fiber $\Phi_B$ consists of all the twisted cubics $A$ that
meet
$B$ in a length-$n$ subscheme $W$.  So $\Phi_B$ has dimension at most
$16-2n$ for $n\le4$, because the various $W$ form an $n$-fold, the
various hyperplanes $G$ containing a fixed $W$ form a $\IP^{4-n}$, and
the various $A$ in a fixed $G$ containing a fixed $W$ form a space of
dimension at most $12-2n$ by Lemma~\Cs4) below.  Hence $R_{3,b,n}$ has
dimension at most $5b+17-2n$, or $n-1$ less than the asserted bound of
$5b+16-n$ for $n\le4$.  Since $R_{3,b,4}$ contains $R_{3,b,n}$ for
$n\ge4$, the latter has dimension at most $5b+9$, which is less than
the asserted bound also for $n=5,6$.  Thus Lemma~\Cs2) holds when $a=3$.

Suppose $a=4$.  Then $b=4$ or $b=5$.  So $A$ is a rational normal curve
in $\IP^4$, and $B$ spans $\IP^4$ too.  Hence $n\le 2b$, because $A$
lies in a quadric hypersurface $Q$ not containing $B$, and $Q\cap B$
has length $2b$.  Now, the fiber $\Phi_B$ consists of all the $A$ that
meet $B$ in a length-$n$ subscheme $W$.  So $\Phi_B$ has dimension at
most $21-2n$ for $n\le6$, because the various $W$ form an $n$-fold, and
the various $A$ containing a fixed $W$ form a space of dimension at
most $21-3n$ by Lemma~\Cs4) below.  Hence $R_{4,b,n}$ has dimension at
most $22+5b-2n$, or $n-1$ less than the asserted bound of $21+5b-n$,
for $n\le6$.  Since the same reasoning as for $R_{4,b,6}$ applies to
$R_{4,b,n}$ for $n\ge4$, the latter has dimension at most $10+5b$,
which is at most the asserted bound also for $7\le n\le 2b$.  Thus
Lemma~\Cs2) holds when $a=4$.  The first assertion of Lemma~\Cs2) is
now proved

To complete the proof of Lemma~\Cs2), consider a pair $(A,B)$ in
$S_{a,n}$.  Then $A$ is a rational normal curve of degree $a$ where
$1\le a\le4$, and $B$ is a six-nodal plane quintic, say with plane $H$.
Then $H\not\supset A$, because there is a quintic hypersurface $F$ that
contains $A$ and $B$, but not $H$, so meets $H$ in $B$.  Fix $B$
momentarily, and let $A$ vary in the fiber $\Phi_B$ of $S_{a,n}$ over
$B$.  First, suppose $a=1$, or $A$ is $n$-secant line of $B$.  Now,
$A\not\subset H$, so $n=1$.  Hence $\dim\Phi_B=4$ since the lines
through a fixed point of $B$ form a $\IP^3$.  Therefore $S_{1,1}$ has
dimension 24, which is the asserted bound.  Thus Lemma~\Cs2) holds when
$a=1$.

Suppose $a=2$.  Then $A$ is an $n$-secant conic of $B$.  Denote the
plane of $A$ by $G$.  Then $G\neq H$ because $H\not\supset A$.  If
$G\cap H$ is a line $L$, then $A\cap B$ is contained in $A\cap L$,
which has length 2; hence, $n\le2$.  If $G\cap H$ is a point, then
clearly $n=1$.  If $n=1$, then $\dim\Phi_B=9$ since, given $P\in B$,
there is a $\IP^4$ of planes $G$ through $P$, and there is a $\IP^4$ of
conics $A$ that lie in a fixed $G$ and pass through $P$.  Hence
$S_{2,1}$ has dimension 29, which is the asserted bound.  If $n=2$,
then $\dim\Phi_B=7$ since the length-2 subschemes $W$ of $B$ form a
surface, the planes $G$ through a fixed $W$ form a $\IP^2$, and
there is a $\IP^3$ of conics $A$ that lie in a fixed $G$ and pass
through $W$.  Hence $S_{2,2}$ has dimension 27, which is one less than
the asserted bound.  Thus Lemma~\Cs2) holds when $a=2$.

Suppose $a=3$.  Then $A$ is a twisted cubic, $n$-secant to $B$.  Denote
the hyperplane of $A$ by $G$.  If $G\supset H$, then $A\cap B$ is
contained in $A\cap H$, which has length 3; hence $n\le3$.  If $G\cap
H$ is a line $L$, then $A\cap L$ has length 2; hence, $n\le2$.  Now,
the fiber $\Phi_B$ has dimension at most $16-2n$ for $n\le3$, because
the length-$n$ subschemes $W$ of $B$ form an $n$-fold, the various
hyperplanes $G$ containing a fixed $W$ form a $\IP^{4-n}$, and the
various $A$ in a fixed $G$ containing a fixed $W$ form a space of
dimension at most $12-2n$ by Lemma~\Cs4) below.  Hence $S_{3,n}$ has
dimension $36-2n$, which is at most the asserted bound for $n=1,2,3$.
Thus Lemma~\Cs2) holds when $a=3$.

Suppose $a=4$.  Then $A$ is a rational normal quartic, $n$-secant to
$B$.  Moreover, $n\le3$; indeed, the plane $H$ of $B$ can intersect $A$
in a subscheme of length at most 3 since $H$ and a general point of $B$
span a hyperplane, which meets $A$ in a subscheme of length 4.  Now,
the fiber $\Phi_B$ has dimension at most $21-2n$, because the
length-$n$ subschemes $W$ of $B$ form an $n$-fold, and the various $A$
containing a fixed $W$ form a space of dimension at most $21-3n$ by
Lemma~\Cs4) below.  Hence $S_{4,n}$ has dimension $41-2n$, which is at
most the asserted bound for $n=1,2,3$.  Thus Lemma~\Cs2) holds when
$a=4$.  The proof of Lemma~\Cs2) is complete, given the following lemma.

\lem4 In the space of all  rational normal curves of degree $a$ in
$\IP^a$, those that contain a given subscheme of length $n$ form a
subset of codimension  $(a-1)n$ for $n\le a+2$.
 \pf
 Consider the incidence scheme $I$ of pairs $(C,W)$ where $C$ is a
rational normal curve and $W$ is a length-$n$ subscheme of $C$.  The
fiber of $I$ over a fixed $C$ is isomorphic to $\IP^n$, and the various
$C$ form an irreducible scheme of dimension $(a+1)^2-4$; hence, $I$ is
irreducible of dimension $(a+1)^2-4+n$.  Consider the projection of $I$
onto the space of $W$, which is a smooth open subscheme $U$ of
dimension $an$ in the Hilbert scheme of $\IP^a$.  This projection is
equivariant for the natural action of $GL(a+1)$.  Since $n\le a+2$,
those pairs $(C,W)$ such that $W$ is concentrated at a point $P$ form a
closed orbit; it is isomorphic to the incidence space of pairs $(C,P)$.
This orbit maps onto a closed orbit in $U$, which has dimension
$(a-1)n+1$ and lies in the closure of every orbit of $U$.  Hence the
fiber over such a $W$ has codimension $(a-1)n$.  Therefore, the fiber
over an arbitrary $W$ has codimension $(a-1)n$ by lower semi-continuity
of codimension.  Thus Lemma~\Cs4) is proved.

 \medbreak

To prove Lemma~\Cs3), first suppose that $B$ is smooth.  Say $B$ spans
an $r$-plane.  Suppose $a=1$.  Then $b\le8$.  Hence $C$ is
(3+b-r)-regular by \cite{GLP, Rmk.~(1), \p.497}, and so 6-regular if
$b\le7$.  Suppose $b=8$.  Then the proof of \cite{GLP, Rmk.~(1)} yields
even that $C$ is 5-regular.  Indeed, by \cite{GLP, Prp.~(1.2), \p.494},
it suffices to show that, on the normalization of $C$, which is the
disjoint union of $A$ and $B$, there is an invertible sheaf $\A$ such
that $h^0(\A)$ is 5 and $H^1(\wedge^2\M\ox\A)$ vanishes, where $\M$ is
the pullback of $\Omega^1_{\IP^4}(1)$.  By the definition of $R_{a,b,n}$,
	$$\Omega^1_{\IP^4}(1)|A\cong\O_{\IP^1}(-1)\oplus\O_{\IP^1}^{3}
	\and \Omega^1_{\IP^4}(1)|B\cong\O_{\IP^1}(-2)^{4}.$$
 So, if we define $\A$ as the structure sheaf on $A$ and as
$\O_{\IP^1}(3)$ on $B$, then $\A$ has the desired two properties.  Thus
Lemma~\Cs3) holds when $a=1$.

Suppose that $a\ge2$.  Then $b\le7$.  We may assume that $a\le b$.
Then $a\le4$.  So $A$ is a rational normal curve.  Hence $C$ is
(4+b-r)-regular by \cite{GLP, Rmk.~(1), \p.497}, and so 6-regular if
$b\le6$.  Suppose $b=7$, so $a=2$.  Then, much as above, the proof of
\cite{GLP, Rmk.~(1)} yields that $C$ is 6-regular.  Indeed,
	$$\displaylines{
   \Omega^1_{\IP^4}(1)|A\cong\O_{\IP^1}(-1)^2\oplus\O_{\IP^1}^{2},\cr
 \Omega^1_{\IP^4}(1)|B\cong\O_{\IP^1}(-2)^{3}\oplus\O_{\IP^1}(-1).\cr}$$
 So, if we define $\A$ as $\O_{\IP^1}(1)$ on $A$ and as $\O_{\IP^1}(3)$
on $B$, then $\A$ has the desired two properties.  Thus Lemma~\Cs3)
holds when $B$ too is smooth.

Suppose that $B$ is a six-nodal plane quintic.  Let $K$ be its plane.
By hypothesis, there is a quintic threefold $F$ that contains both $A$
and $B$, but not $K$.   Hence  $F\cap K=B$ by Bezout's theorem.  Therefore,
	$$A\cap K\subset A\cap F\cap K\subset B.$$
 Now, $A$ is a rational normal curve of degree $a\le4$.  Hence, $A$ is
smooth and 2-regular, so 4-regular.  Therefore, $C$ is 5-regular by the
following lemma, and after it is proved, the proofs of Lemma~\Cs3) and
of Theorem~\Cs1) will be complete.

\lem5 Let $A$ and $B$ be two reduced curves in $\IP^4$, and $C$ their
union.  Assume that $A$ and $B$ have no common components, that $B$
spans a plane $K$, and that $A$ is smooth along $A\cap K$.  If $A$ is
$(m-1)$-regular where $m\ge\deg B$, then $C$ is $m$-regular.
 \pf
 First, we'll use Hirschowitz's ``m\'ethode d'Horace'' (see
\cite{Hirsch}).  Let $H$ be a general hyperplane containing $K$, and
set $D:=C\cap H$.  Then there is a natural exact sequence of sheaves,
	$$0\to \I_{A/\IP^4}(-1)\to\I_{C/\IP^4}\TO u\I_{D/H}\to0$$
  because
	$$\eqalign{\Ker u&=\I_{C/\IP^4}\cap\I_{H/\IP^4}
	=\I_{A/\IP^4}\cap\I_{B/\IP^4}\cap\I_{H/\IP^4}\cr
	&=\I_{A/\IP^4}\cap\I_{H/\IP^4}=\I_{A/\IP^4}(-1);\cr}$$
 the last equation holds because $A$ is reduced and because $H$, being
general, contains no component of $A$.  By hypothesis, $A$ is
$(m-1)$-regular.  Therefore, $C$ will be $m$-regular if $\I_{D/H}$ is
so.

Set $A':=A\cap H$.  Then $A'$ is the disjoint union of two finite
closed subschemes: the part $S$ of $A'$ with support off $K$, and the
part $T$ with support on $K$.  Correspondingly, there is a natural
exact sequence of sheaves,
	$$0\to \I_{A'/H}\to\I_{S/H}\oplus\I_{T/H}\to\O_H\to0$$
 By hypothesis, $\I_{A/\IP^4}$ is $(m-1)$-regular; whence, so is
$\I_{A'/H}$.  On the other hand, $\O_H$ is 0-regular, so
$(m-1)$-regular.  Therefore, $\I_{S/H}$ is $(m-1)$-regular.

Since $A$ is smooth along $A\cap K$, the latter is a locally principal
subscheme of $A$.  Hence, since $H$ is general containing $K$, the
schemes $A\cap K$ and $T$ coincide.  By hypothesis, $A\cap K$ is a
subscheme of $B$.  Hence is $T$ also.  It follows that the scheme $D$
is the disjoint union of $S$ and $B$.  Indeed, this statement is
clearly true off $T$.  So work locally at an arbitrary point $P$ of
$T$.  Since $D\supset B\supset T$, we have
	$$(\I_{C/\IP^4,P}+\I_{H/\IP^4,P})\subset\I_{B/\IP^4,P}
		\subset(\I_{A/\IP^4,P}+\I_{H/\IP^4,P}),\eqlt51$$
 and it suffices to check that the first inclusion is an equality.
Given $\beta\in\I_{B/\IP^4,P}$, write $\beta=\alpha+\gamma$ where
$\alpha\in\I_{A/\IP^4,P}$ and $\gamma\in\I_{H/\IP^4,P}$.  Then
	$$\alpha=\beta-\gamma\in\I_{A/\IP^4,P}\cap\I_{B/\IP^4,P}
	=\I_{C/\IP^4,P}.$$
 Hence $\beta$ lies in the left hand term of \Cd51), as required.  Thus
$D=S\cup B$.  Therefore, $D\cap K=B$ because $K$ is disjoint from $S$
and contains $B$.

Again we'll use Hirschowitz's ``m\'ethode d'Horace.''  Much as above,
there is a natural exact sequence of sheaves,
	$$0\to \I_{S/H}(-1)\to\I_{D/H}\TO v\I_{B/K}\to0,$$
  because $D\cap K=B$ and because
     $$\eqalign{\Ker v&=\I_{D/H}\cap\I_{K/H}
     =\I_{S/H}\cap\I_{B/H}\cap\I_{K/H}\cr
     &=\I_{S/H}\cap\I_{K/H}=\I_{S/H}(-1);\cr}$$
 the last equation holds because $K$ is disjoint from $S$.  Now, it was
proved above that $\I_{S/H}$ is $(m-1)$-regular.  On the other hand,
$\I_{B/K}$ is equal to $\O_K(-b)$ where $b:=\deg B$, and the latter
sheaf is $m$-regular since $m\ge b$.  Therefore, $\I_{D/H}$ is $m$-regular.
Consequently, $C$ is $m$-regular by the first paragraph, and the proof
is complete.

\newsect References

\references

ACGH
E. Arbarello, M. Cornalba, P. Griffiths, and J. Harris
 \book Algebraic Curves \bookinfo Grundlehren {\bf 267}, Springer-Verlag, 1985

BE
E. Ballico and P. Ellia
 \paper On the Postulation of Curves in $\IP^4$
 \mathz 188 1985 215--23

Cetal
 P. Candelas, P.S.Green, X.C. de la Ossa, and L. Parkes
 \paper A pair of Calabi-Yau manifolds as an exactly soluble
superconformal theory
 \np 358 1991 21--74

Cil
C. Ciliberto
 \paper Hilbert functions of finite sets of points and the genus of
a curve in a projective space
 \RdiPapa 24--73

C
H. Clemens
 \paper Homological equivalece, modulo algebraic equivalence, is not
fin\-ite\-ly generated
 \pmihes 58 1983 19--38

C86
H. Clemens
 \paper Curves on higher-dimensional complex projective manifolds
  \paperinfo Proc. International Cong. Math., Berkeley, 1986, 634--40

CJ
M. Coppens and T. Johnsen
 \paper Secant lines of smooth projective curves
  \Zeuthen 61--87

dA
J. d'Almeida
 \paper Courbes de l'espace projectif: S\'eries lin\'eaires
incompl\`etes et multis\'ecantes
 \crelle 370 1986 30--51

ES
 G. Ellingsrud and S.A. Str\o mme
 \paper The number of twisted cubic curves on the general
quintic threefold
 \ms 76 1995 5--34

GP
L. Gruson and C. Peskine
 \paper Genre des courbes de l'espace projectif
 \Tromso 31-60

GLP
L. Gruson, R. Lazarsfeld, and C. Peskine
 \paper On a theorem of Castelnuovo and the equations defining space
curves
 \invent 72 1983 491--506

Gu
M. Guest
 \paper On the space of holomorphic maps from the Riemann sphere to the
quadric cone.
 \qjm 45 1994 57--75

H
R. Hartshorne
 \book Algebraic Geometry
 \bookinfo GTM {\bf 52}, Springer-Verlag, 1977

Hirsch
A. Hirschowitz
 \paper La m\'ethode d'Horace pour l\'interpolation \'a plusieurs
variables
 \mm 50 1985 337--88

Huyb
D. Huybrechts
 \paper A note on the finiteness of smooth rational curves on
Calabi--Yau threefolds
 \paperinfo Preprint, Max-Planck-Institut, April 14, 1993

Katz
S. Katz
 \paper On the finiteness of rational curves on quintic threefolds
 \comp 60 1986 151--62

Kont
M. Kontsevich
 \paper Enumeration of rational curves via torus actions
 \preprint hep-th/9402147

L
A. Laudal
\paper A generalized trisecant lemma
\Tromso 112--49

LY
B. Lian and S.-T. Yau
 \paper Arithmetic Properties of Mirror Map and Quantum Coupling
 \preprint hep-th/9411234

LT
 A. Libgober and J.Teitelbaum
 \paper Lines on Calabi-Yau complete intersections,
mirror symmetry and Picard-Fuchs equations
 \duke 69(1) 1993 29--39

N
P. Nijsse
 \paper Clemens' conjecture for octic and nonic curves
 \paperinfo Preprint, U. Leiden, 1993, to appear in Indag. Math

Oguiso
K. Oguiso
 \paper Two remarks on Calabai--Yau Moishezon threefolds
 \preprint MPI 93-15

Ramella
L. Ramella
 \paper La stratification du sch\'ema de Hilbert des courbes rationelles
de $\IP^n$ par le fibr\'e tangent restreint
 \crasp 311 1990 181--4

Rath
J. Rathmann
 \paper The uniform position principle for curves in characteristic p
 \ma 276 1987 565--79

V
I. Vainsencher
 \paper Enumeration of  $n$-fold tangent hyperplanes to a surface
 \jag 4 1995 503--26

Verdier
 J.-L. Verdier
 \paper Two dimensional sigma-models and harmonic maps from ${\bf S}^2$
to ${\bf S}^{2n}$
 \paperinfo in ``Group Theoretical Methods in Physics $\bullet$
Proceedings, Istanbul, Turkey,'' Springer Lecture Notes in Physics {\bf
180} (1983), 136--41

Westhoff
R. Westhoff
 \paper Curves and normal functions on higher dimensional complex
varieties \paperinfo PhD Thesis, U. Utah, June 1993

\endreferences

  \bigskip

 \eightpoint\smc Mathematical Institute, University of Bergen, All\'egaten
55, N-5007 Berg\-en, Norway
 \medskip
 Department of Mathematics, 2--278 MIT, Cambridge, MA 02139,
U.S.A.

\bye